\documentclass[10pt]{article}
\usepackage{pdfsync}
\usepackage{amsmath}
\usepackage{amsthm}

\usepackage{authblk}

\usepackage{graphicx}
\usepackage{color}
\definecolor{AliceBlue}{rgb}{0.94,0.97,1.00}
\definecolor{AntiqueWhite1}{rgb}{1.00,0.94,0.86}
\definecolor{AntiqueWhite2}{rgb}{0.93,0.87,0.80}
\definecolor{AntiqueWhite3}{rgb}{0.80,0.75,0.69}
\definecolor{AntiqueWhite4}{rgb}{0.55,0.51,0.47}
\definecolor{AntiqueWhite}{rgb}{0.98,0.92,0.84}
\definecolor{BlanchedAlmond}{rgb}{1.00,0.92,0.80}
\definecolor{BlueViolet}{rgb}{0.54,0.17,0.89}
\definecolor{CadetBlue1}{rgb}{0.60,0.96,1.00}
\definecolor{CadetBlue2}{rgb}{0.56,0.90,0.93}
\definecolor{CadetBlue3}{rgb}{0.48,0.77,0.80}
\definecolor{CadetBlue4}{rgb}{0.33,0.53,0.55}
\definecolor{CadetBlue}{rgb}{0.37,0.62,0.63}
\definecolor{CornflowerBlue}{rgb}{0.39,0.58,0.93}
\definecolor{DarkBlue}{rgb}{0.00,0.00,0.55}
\definecolor{DarkCyan}{rgb}{0.00,0.55,0.55}
\definecolor{DarkGoldenrod1}{rgb}{1.00,0.73,0.06}
\definecolor{DarkGoldenrod2}{rgb}{0.93,0.68,0.05}
\definecolor{DarkGoldenrod3}{rgb}{0.80,0.58,0.05}
\definecolor{DarkGoldenrod4}{rgb}{0.55,0.40,0.03}
\definecolor{DarkGoldenrod}{rgb}{0.72,0.53,0.04}
\definecolor{DarkGray}{rgb}{0.66,0.66,0.66}
\definecolor{DarkGreen}{rgb}{0.00,0.39,0.00}
\definecolor{DarkGrey}{rgb}{0.66,0.66,0.66}
\definecolor{DarkKhaki}{rgb}{0.74,0.72,0.42}
\definecolor{DarkMagenta}{rgb}{0.55,0.00,0.55}
\definecolor{DarkOliveGreen1}{rgb}{0.79,1.00,0.44}
\definecolor{DarkOliveGreen2}{rgb}{0.74,0.93,0.41}
\definecolor{DarkOliveGreen3}{rgb}{0.64,0.80,0.35}
\definecolor{DarkOliveGreen4}{rgb}{0.43,0.55,0.24}
\definecolor{DarkOliveGreen}{rgb}{0.33,0.42,0.18}
\definecolor{DarkOrange1}{rgb}{1.00,0.50,0.00}
\definecolor{DarkOrange2}{rgb}{0.93,0.46,0.00}
\definecolor{DarkOrange3}{rgb}{0.80,0.40,0.00}
\definecolor{DarkOrange4}{rgb}{0.55,0.27,0.00}
\definecolor{DarkOrange}{rgb}{1.00,0.55,0.00}
\definecolor{DarkOrchid1}{rgb}{0.75,0.24,1.00}
\definecolor{DarkOrchid2}{rgb}{0.70,0.23,0.93}
\definecolor{DarkOrchid3}{rgb}{0.60,0.20,0.80}
\definecolor{DarkOrchid4}{rgb}{0.41,0.13,0.55}
\definecolor{DarkOrchid}{rgb}{0.60,0.20,0.80}
\definecolor{DarkRed}{rgb}{0.55,0.00,0.00}
\definecolor{DarkSalmon}{rgb}{0.91,0.59,0.48}
\definecolor{DarkSeaGreen1}{rgb}{0.76,1.00,0.76}
\definecolor{DarkSeaGreen2}{rgb}{0.71,0.93,0.71}
\definecolor{DarkSeaGreen3}{rgb}{0.61,0.80,0.61}
\definecolor{DarkSeaGreen4}{rgb}{0.41,0.55,0.41}
\definecolor{DarkSeaGreen}{rgb}{0.56,0.74,0.56}
\definecolor{DarkSlateBlue}{rgb}{0.28,0.24,0.55}
\definecolor{DarkSlateGray1}{rgb}{0.59,1.00,1.00}
\definecolor{DarkSlateGray2}{rgb}{0.55,0.93,0.93}
\definecolor{DarkSlateGray3}{rgb}{0.47,0.80,0.80}
\definecolor{DarkSlateGray4}{rgb}{0.32,0.55,0.55}
\definecolor{DarkSlateGray}{rgb}{0.18,0.31,0.31}
\definecolor{DarkSlateGrey}{rgb}{0.18,0.31,0.31}
\definecolor{DarkTurquoise}{rgb}{0.00,0.81,0.82}
\definecolor{DarkViolet}{rgb}{0.58,0.00,0.83}
\definecolor{DeepPink1}{rgb}{1.00,0.08,0.58}
\definecolor{DeepPink2}{rgb}{0.93,0.07,0.54}
\definecolor{DeepPink3}{rgb}{0.80,0.06,0.46}
\definecolor{DeepPink4}{rgb}{0.55,0.04,0.31}
\definecolor{DeepPink}{rgb}{1.00,0.08,0.58}
\definecolor{DeepSkyBlue1}{rgb}{0.00,0.75,1.00}
\definecolor{DeepSkyBlue2}{rgb}{0.00,0.70,0.93}
\definecolor{DeepSkyBlue3}{rgb}{0.00,0.60,0.80}
\definecolor{DeepSkyBlue4}{rgb}{0.00,0.41,0.55}
\definecolor{DeepSkyBlue}{rgb}{0.00,0.75,1.00}
\definecolor{DimGray}{rgb}{0.41,0.41,0.41}
\definecolor{DimGrey}{rgb}{0.41,0.41,0.41}
\definecolor{DodgerBlue1}{rgb}{0.12,0.56,1.00}
\definecolor{DodgerBlue2}{rgb}{0.11,0.53,0.93}
\definecolor{DodgerBlue3}{rgb}{0.09,0.45,0.80}
\definecolor{DodgerBlue4}{rgb}{0.06,0.31,0.55}
\definecolor{DodgerBlue}{rgb}{0.12,0.56,1.00}
\definecolor{FloralWhite}{rgb}{1.00,0.98,0.94}
\definecolor{ForestGreen}{rgb}{0.13,0.55,0.13}
\definecolor{GhostWhite}{rgb}{0.97,0.97,1.00}
\definecolor{GreenYellow}{rgb}{0.68,1.00,0.18}
\definecolor{HotPink1}{rgb}{1.00,0.43,0.71}
\definecolor{HotPink2}{rgb}{0.93,0.42,0.65}
\definecolor{HotPink3}{rgb}{0.80,0.38,0.56}
\definecolor{HotPink4}{rgb}{0.55,0.23,0.38}
\definecolor{HotPink}{rgb}{1.00,0.41,0.71}
\definecolor{IndianRed1}{rgb}{1.00,0.42,0.42}
\definecolor{IndianRed2}{rgb}{0.93,0.39,0.39}
\definecolor{IndianRed3}{rgb}{0.80,0.33,0.33}
\definecolor{IndianRed4}{rgb}{0.55,0.23,0.23}
\definecolor{IndianRed}{rgb}{0.80,0.36,0.36}
\definecolor{LavenderBlush1}{rgb}{1.00,0.94,0.96}
\definecolor{LavenderBlush2}{rgb}{0.93,0.88,0.90}
\definecolor{LavenderBlush3}{rgb}{0.80,0.76,0.77}
\definecolor{LavenderBlush4}{rgb}{0.55,0.51,0.53}
\definecolor{LavenderBlush}{rgb}{1.00,0.94,0.96}
\definecolor{LawnGreen}{rgb}{0.49,0.99,0.00}
\definecolor{LemonChiffon1}{rgb}{1.00,0.98,0.80}
\definecolor{LemonChiffon2}{rgb}{0.93,0.91,0.75}
\definecolor{LemonChiffon3}{rgb}{0.80,0.79,0.65}
\definecolor{LemonChiffon4}{rgb}{0.55,0.54,0.44}
\definecolor{LemonChiffon}{rgb}{1.00,0.98,0.80}
\definecolor{LightBlue1}{rgb}{0.75,0.94,1.00}
\definecolor{LightBlue2}{rgb}{0.70,0.87,0.93}
\definecolor{LightBlue3}{rgb}{0.60,0.75,0.80}
\definecolor{LightBlue4}{rgb}{0.41,0.51,0.55}
\definecolor{LightBlue}{rgb}{0.68,0.85,0.90}
\definecolor{LightCoral}{rgb}{0.94,0.50,0.50}
\definecolor{LightCyan1}{rgb}{0.88,1.00,1.00}
\definecolor{LightCyan2}{rgb}{0.82,0.93,0.93}
\definecolor{LightCyan3}{rgb}{0.71,0.80,0.80}
\definecolor{LightCyan4}{rgb}{0.48,0.55,0.55}
\definecolor{LightCyan}{rgb}{0.88,1.00,1.00}
\definecolor{LightGoldenrod1}{rgb}{1.00,0.93,0.55}
\definecolor{LightGoldenrod2}{rgb}{0.93,0.86,0.51}
\definecolor{LightGoldenrod3}{rgb}{0.80,0.75,0.44}
\definecolor{LightGoldenrod4}{rgb}{0.55,0.51,0.30}
\definecolor{LightGoldenrodYellow}{rgb}{0.98,0.98,0.82}
\definecolor{LightGoldenrod}{rgb}{0.93,0.87,0.51}
\definecolor{LightGray}{rgb}{0.83,0.83,0.83}
\definecolor{LightGreen}{rgb}{0.56,0.93,0.56}
\definecolor{LightGrey}{rgb}{0.83,0.83,0.83}
\definecolor{LightPink1}{rgb}{1.00,0.68,0.73}
\definecolor{LightPink2}{rgb}{0.93,0.64,0.68}
\definecolor{LightPink3}{rgb}{0.80,0.55,0.58}
\definecolor{LightPink4}{rgb}{0.55,0.37,0.40}
\definecolor{LightPink}{rgb}{1.00,0.71,0.76}
\definecolor{LightSalmon1}{rgb}{1.00,0.63,0.48}
\definecolor{LightSalmon2}{rgb}{0.93,0.58,0.45}
\definecolor{LightSalmon3}{rgb}{0.80,0.51,0.38}
\definecolor{LightSalmon4}{rgb}{0.55,0.34,0.26}
\definecolor{LightSalmon}{rgb}{1.00,0.63,0.48}
\definecolor{LightSeaGreen}{rgb}{0.13,0.70,0.67}
\definecolor{LightSkyBlue1}{rgb}{0.69,0.89,1.00}
\definecolor{LightSkyBlue2}{rgb}{0.64,0.83,0.93}
\definecolor{LightSkyBlue3}{rgb}{0.55,0.71,0.80}
\definecolor{LightSkyBlue4}{rgb}{0.38,0.48,0.55}
\definecolor{LightSkyBlue}{rgb}{0.53,0.81,0.98}
\definecolor{LightSlateBlue}{rgb}{0.52,0.44,1.00}
\definecolor{LightSlateGray}{rgb}{0.47,0.53,0.60}
\definecolor{LightSlateGrey}{rgb}{0.47,0.53,0.60}
\definecolor{LightSteelBlue1}{rgb}{0.79,0.88,1.00}
\definecolor{LightSteelBlue2}{rgb}{0.74,0.82,0.93}
\definecolor{LightSteelBlue3}{rgb}{0.64,0.71,0.80}
\definecolor{LightSteelBlue4}{rgb}{0.43,0.48,0.55}
\definecolor{LightSteelBlue}{rgb}{0.69,0.77,0.87}
\definecolor{LightYellow1}{rgb}{1.00,1.00,0.88}
\definecolor{LightYellow2}{rgb}{0.93,0.93,0.82}
\definecolor{LightYellow3}{rgb}{0.80,0.80,0.71}
\definecolor{LightYellow4}{rgb}{0.55,0.55,0.48}
\definecolor{LightYellow}{rgb}{1.00,1.00,0.88}
\definecolor{LimeGreen}{rgb}{0.20,0.80,0.20}
\definecolor{MediumAquamarine}{rgb}{0.40,0.80,0.67}
\definecolor{MediumBlue}{rgb}{0.00,0.00,0.80}
\definecolor{MediumOrchid1}{rgb}{0.88,0.40,1.00}
\definecolor{MediumOrchid2}{rgb}{0.82,0.37,0.93}
\definecolor{MediumOrchid3}{rgb}{0.71,0.32,0.80}
\definecolor{MediumOrchid4}{rgb}{0.48,0.22,0.55}
\definecolor{MediumOrchid}{rgb}{0.73,0.33,0.83}
\definecolor{MediumPurple1}{rgb}{0.67,0.51,1.00}
\definecolor{MediumPurple2}{rgb}{0.62,0.47,0.93}
\definecolor{MediumPurple3}{rgb}{0.54,0.41,0.80}
\definecolor{MediumPurple4}{rgb}{0.36,0.28,0.55}
\definecolor{MediumPurple}{rgb}{0.58,0.44,0.86}
\definecolor{MediumSeaGreen}{rgb}{0.24,0.70,0.44}
\definecolor{MediumSlateBlue}{rgb}{0.48,0.41,0.93}
\definecolor{MediumSpringGreen}{rgb}{0.00,0.98,0.60}
\definecolor{MediumTurquoise}{rgb}{0.28,0.82,0.80}
\definecolor{MediumVioletRed}{rgb}{0.78,0.08,0.52}
\definecolor{MidnightBlue}{rgb}{0.10,0.10,0.44}
\definecolor{MintCream}{rgb}{0.96,1.00,0.98}
\definecolor{MistyRose1}{rgb}{1.00,0.89,0.88}
\definecolor{MistyRose2}{rgb}{0.93,0.84,0.82}
\definecolor{MistyRose3}{rgb}{0.80,0.72,0.71}
\definecolor{MistyRose4}{rgb}{0.55,0.49,0.48}
\definecolor{MistyRose}{rgb}{1.00,0.89,0.88}
\definecolor{NavajoWhite1}{rgb}{1.00,0.87,0.68}
\definecolor{NavajoWhite2}{rgb}{0.93,0.81,0.63}
\definecolor{NavajoWhite3}{rgb}{0.80,0.70,0.55}
\definecolor{NavajoWhite4}{rgb}{0.55,0.47,0.37}
\definecolor{NavajoWhite}{rgb}{1.00,0.87,0.68}
\definecolor{NavyBlue}{rgb}{0.00,0.00,0.50}
\definecolor{OldLace}{rgb}{0.99,0.96,0.90}
\definecolor{OliveDrab1}{rgb}{0.75,1.00,0.24}
\definecolor{OliveDrab2}{rgb}{0.70,0.93,0.23}
\definecolor{OliveDrab3}{rgb}{0.60,0.80,0.20}
\definecolor{OliveDrab4}{rgb}{0.41,0.55,0.13}
\definecolor{OliveDrab}{rgb}{0.42,0.56,0.14}
\definecolor{OrangeRed1}{rgb}{1.00,0.27,0.00}
\definecolor{OrangeRed2}{rgb}{0.93,0.25,0.00}
\definecolor{OrangeRed3}{rgb}{0.80,0.22,0.00}
\definecolor{OrangeRed4}{rgb}{0.55,0.15,0.00}
\definecolor{OrangeRed}{rgb}{1.00,0.27,0.00}
\definecolor{PaleGoldenrod}{rgb}{0.93,0.91,0.67}
\definecolor{PaleGreen1}{rgb}{0.60,1.00,0.60}
\definecolor{PaleGreen2}{rgb}{0.56,0.93,0.56}
\definecolor{PaleGreen3}{rgb}{0.49,0.80,0.49}
\definecolor{PaleGreen4}{rgb}{0.33,0.55,0.33}
\definecolor{PaleGreen}{rgb}{0.60,0.98,0.60}
\definecolor{PaleTurquoise1}{rgb}{0.73,1.00,1.00}
\definecolor{PaleTurquoise2}{rgb}{0.68,0.93,0.93}
\definecolor{PaleTurquoise3}{rgb}{0.59,0.80,0.80}
\definecolor{PaleTurquoise4}{rgb}{0.40,0.55,0.55}
\definecolor{PaleTurquoise}{rgb}{0.69,0.93,0.93}
\definecolor{PaleVioletRed1}{rgb}{1.00,0.51,0.67}
\definecolor{PaleVioletRed2}{rgb}{0.93,0.47,0.62}
\definecolor{PaleVioletRed3}{rgb}{0.80,0.41,0.54}
\definecolor{PaleVioletRed4}{rgb}{0.55,0.28,0.36}
\definecolor{PaleVioletRed}{rgb}{0.86,0.44,0.58}
\definecolor{PapayaWhip}{rgb}{1.00,0.94,0.84}
\definecolor{PeachPuff1}{rgb}{1.00,0.85,0.73}
\definecolor{PeachPuff2}{rgb}{0.93,0.80,0.68}
\definecolor{PeachPuff3}{rgb}{0.80,0.69,0.58}
\definecolor{PeachPuff4}{rgb}{0.55,0.47,0.40}
\definecolor{PeachPuff}{rgb}{1.00,0.85,0.73}
\definecolor{PowderBlue}{rgb}{0.69,0.88,0.90}
\definecolor{RosyBrown1}{rgb}{1.00,0.76,0.76}
\definecolor{RosyBrown2}{rgb}{0.93,0.71,0.71}
\definecolor{RosyBrown3}{rgb}{0.80,0.61,0.61}
\definecolor{RosyBrown4}{rgb}{0.55,0.41,0.41}
\definecolor{RosyBrown}{rgb}{0.74,0.56,0.56}
\definecolor{RoyalBlue1}{rgb}{0.28,0.46,1.00}
\definecolor{RoyalBlue2}{rgb}{0.26,0.43,0.93}
\definecolor{RoyalBlue3}{rgb}{0.23,0.37,0.80}
\definecolor{RoyalBlue4}{rgb}{0.15,0.25,0.55}
\definecolor{RoyalBlue}{rgb}{0.25,0.41,0.88}
\definecolor{SaddleBrown}{rgb}{0.55,0.27,0.07}
\definecolor{SandyBrown}{rgb}{0.96,0.64,0.38}
\definecolor{SeaGreen1}{rgb}{0.33,1.00,0.62}
\definecolor{SeaGreen2}{rgb}{0.31,0.93,0.58}
\definecolor{SeaGreen3}{rgb}{0.26,0.80,0.50}
\definecolor{SeaGreen4}{rgb}{0.18,0.55,0.34}
\definecolor{SeaGreen}{rgb}{0.18,0.55,0.34}
\definecolor{SkyBlue1}{rgb}{0.53,0.81,1.00}
\definecolor{SkyBlue2}{rgb}{0.49,0.75,0.93}
\definecolor{SkyBlue3}{rgb}{0.42,0.65,0.80}
\definecolor{SkyBlue4}{rgb}{0.29,0.44,0.55}
\definecolor{SkyBlue}{rgb}{0.53,0.81,0.92}
\definecolor{SlateBlue1}{rgb}{0.51,0.44,1.00}
\definecolor{SlateBlue2}{rgb}{0.48,0.40,0.93}
\definecolor{SlateBlue3}{rgb}{0.41,0.35,0.80}
\definecolor{SlateBlue4}{rgb}{0.28,0.24,0.55}
\definecolor{SlateBlue}{rgb}{0.42,0.35,0.80}
\definecolor{SlateGray1}{rgb}{0.78,0.89,1.00}
\definecolor{SlateGray2}{rgb}{0.73,0.83,0.93}
\definecolor{SlateGray3}{rgb}{0.62,0.71,0.80}
\definecolor{SlateGray4}{rgb}{0.42,0.48,0.55}
\definecolor{SlateGray}{rgb}{0.44,0.50,0.56}
\definecolor{SlateGrey}{rgb}{0.44,0.50,0.56}
\definecolor{SpringGreen1}{rgb}{0.00,1.00,0.50}
\definecolor{SpringGreen2}{rgb}{0.00,0.93,0.46}
\definecolor{SpringGreen3}{rgb}{0.00,0.80,0.40}
\definecolor{SpringGreen4}{rgb}{0.00,0.55,0.27}
\definecolor{SpringGreen}{rgb}{0.00,1.00,0.50}
\definecolor{SteelBlue1}{rgb}{0.39,0.72,1.00}
\definecolor{SteelBlue2}{rgb}{0.36,0.67,0.93}
\definecolor{SteelBlue3}{rgb}{0.31,0.58,0.80}
\definecolor{SteelBlue4}{rgb}{0.21,0.39,0.55}
\definecolor{SteelBlue}{rgb}{0.27,0.51,0.71}
\definecolor{VioletRed1}{rgb}{1.00,0.24,0.59}
\definecolor{VioletRed2}{rgb}{0.93,0.23,0.55}
\definecolor{VioletRed3}{rgb}{0.80,0.20,0.47}
\definecolor{VioletRed4}{rgb}{0.55,0.13,0.32}
\definecolor{VioletRed}{rgb}{0.82,0.13,0.56}
\definecolor{WhiteSmoke}{rgb}{0.96,0.96,0.96}
\definecolor{YellowGreen}{rgb}{0.60,0.80,0.20}
\definecolor{aliceblue}{rgb}{0.94,0.97,1.00}
\definecolor{antiquewhite}{rgb}{0.98,0.92,0.84}
\definecolor{aquamarine1}{rgb}{0.50,1.00,0.83}
\definecolor{aquamarine2}{rgb}{0.46,0.93,0.78}
\definecolor{aquamarine3}{rgb}{0.40,0.80,0.67}
\definecolor{aquamarine4}{rgb}{0.27,0.55,0.45}
\definecolor{aquamarine}{rgb}{0.50,1.00,0.83}
\definecolor{azure1}{rgb}{0.94,1.00,1.00}
\definecolor{azure2}{rgb}{0.88,0.93,0.93}
\definecolor{azure3}{rgb}{0.76,0.80,0.80}
\definecolor{azure4}{rgb}{0.51,0.55,0.55}
\definecolor{azure}{rgb}{0.94,1.00,1.00}
\definecolor{beige}{rgb}{0.96,0.96,0.86}
\definecolor{bisque1}{rgb}{1.00,0.89,0.77}
\definecolor{bisque2}{rgb}{0.93,0.84,0.72}
\definecolor{bisque3}{rgb}{0.80,0.72,0.62}
\definecolor{bisque4}{rgb}{0.55,0.49,0.42}
\definecolor{bisque}{rgb}{1.00,0.89,0.77}
\definecolor{black}{rgb}{0.00,0.00,0.00}
\definecolor{blanchedalmond}{rgb}{1.00,0.92,0.80}
\definecolor{blue1}{rgb}{0.00,0.00,1.00}
\definecolor{blue2}{rgb}{0.00,0.00,0.93}
\definecolor{blue3}{rgb}{0.00,0.00,0.80}
\definecolor{blue4}{rgb}{0.00,0.00,0.55}
\definecolor{blueviolet}{rgb}{0.54,0.17,0.89}
\definecolor{blue}{rgb}{0.00,0.00,1.00}
\definecolor{brown1}{rgb}{1.00,0.25,0.25}
\definecolor{brown2}{rgb}{0.93,0.23,0.23}
\definecolor{brown3}{rgb}{0.80,0.20,0.20}
\definecolor{brown4}{rgb}{0.55,0.14,0.14}
\definecolor{brown}{rgb}{0.65,0.16,0.16}
\definecolor{burlywood1}{rgb}{1.00,0.83,0.61}
\definecolor{burlywood2}{rgb}{0.93,0.77,0.57}
\definecolor{burlywood3}{rgb}{0.80,0.67,0.49}
\definecolor{burlywood4}{rgb}{0.55,0.45,0.33}
\definecolor{burlywood}{rgb}{0.87,0.72,0.53}
\definecolor{cadetblue}{rgb}{0.37,0.62,0.63}
\definecolor{chartreuse1}{rgb}{0.50,1.00,0.00}
\definecolor{chartreuse2}{rgb}{0.46,0.93,0.00}
\definecolor{chartreuse3}{rgb}{0.40,0.80,0.00}
\definecolor{chartreuse4}{rgb}{0.27,0.55,0.00}
\definecolor{chartreuse}{rgb}{0.50,1.00,0.00}
\definecolor{chocolate1}{rgb}{1.00,0.50,0.14}
\definecolor{chocolate2}{rgb}{0.93,0.46,0.13}
\definecolor{chocolate3}{rgb}{0.80,0.40,0.11}
\definecolor{chocolate4}{rgb}{0.55,0.27,0.07}
\definecolor{chocolate}{rgb}{0.82,0.41,0.12}
\definecolor{coral1}{rgb}{1.00,0.45,0.34}
\definecolor{coral2}{rgb}{0.93,0.42,0.31}
\definecolor{coral3}{rgb}{0.80,0.36,0.27}
\definecolor{coral4}{rgb}{0.55,0.24,0.18}
\definecolor{coral}{rgb}{1.00,0.50,0.31}
\definecolor{cornflowerblue}{rgb}{0.39,0.58,0.93}
\definecolor{cornsilk1}{rgb}{1.00,0.97,0.86}
\definecolor{cornsilk2}{rgb}{0.93,0.91,0.80}
\definecolor{cornsilk3}{rgb}{0.80,0.78,0.69}
\definecolor{cornsilk4}{rgb}{0.55,0.53,0.47}
\definecolor{cornsilk}{rgb}{1.00,0.97,0.86}
\definecolor{cyan1}{rgb}{0.00,1.00,1.00}
\definecolor{cyan2}{rgb}{0.00,0.93,0.93}
\definecolor{cyan3}{rgb}{0.00,0.80,0.80}
\definecolor{cyan4}{rgb}{0.00,0.55,0.55}
\definecolor{cyan}{rgb}{0.00,1.00,1.00}
\definecolor{darkblue}{rgb}{0.00,0.00,0.55}
\definecolor{darkcyan}{rgb}{0.00,0.55,0.55}
\definecolor{darkgoldenrod}{rgb}{0.72,0.53,0.04}
\definecolor{darkgray}{rgb}{0.66,0.66,0.66}
\definecolor{darkgreen}{rgb}{0.00,0.39,0.00}
\definecolor{darkgrey}{rgb}{0.66,0.66,0.66}
\definecolor{darkkhaki}{rgb}{0.74,0.72,0.42}
\definecolor{darkmagenta}{rgb}{0.55,0.00,0.55}
\definecolor{darkolive}{rgb}{0.33,0.42,0.18}
\definecolor{darkorange}{rgb}{1.00,0.55,0.00}
\definecolor{darkorchid}{rgb}{0.60,0.20,0.80}
\definecolor{darkred}{rgb}{0.55,0.00,0.00}
\definecolor{darksalmon}{rgb}{0.91,0.59,0.48}
\definecolor{darksea}{rgb}{0.56,0.74,0.56}
\definecolor{darkslate}{rgb}{0.18,0.31,0.31}
\definecolor{darkslate}{rgb}{0.18,0.31,0.31}
\definecolor{darkslate}{rgb}{0.28,0.24,0.55}
\definecolor{darkturquoise}{rgb}{0.00,0.81,0.82}
\definecolor{darkviolet}{rgb}{0.58,0.00,0.83}
\definecolor{deeppink}{rgb}{1.00,0.08,0.58}
\definecolor{deepsky}{rgb}{0.00,0.75,1.00}
\definecolor{dimgray}{rgb}{0.41,0.41,0.41}
\definecolor{dimgrey}{rgb}{0.41,0.41,0.41}
\definecolor{dodgerblue}{rgb}{0.12,0.56,1.00}
\definecolor{firebrick1}{rgb}{1.00,0.19,0.19}
\definecolor{firebrick2}{rgb}{0.93,0.17,0.17}
\definecolor{firebrick3}{rgb}{0.80,0.15,0.15}
\definecolor{firebrick4}{rgb}{0.55,0.10,0.10}
\definecolor{firebrick}{rgb}{0.70,0.13,0.13}
\definecolor{floralwhite}{rgb}{1.00,0.98,0.94}
\definecolor{forestgreen}{rgb}{0.13,0.55,0.13}
\definecolor{gainsboro}{rgb}{0.86,0.86,0.86}
\definecolor{ghostwhite}{rgb}{0.97,0.97,1.00}
\definecolor{gold1}{rgb}{1.00,0.84,0.00}
\definecolor{gold2}{rgb}{0.93,0.79,0.00}
\definecolor{gold3}{rgb}{0.80,0.68,0.00}
\definecolor{gold4}{rgb}{0.55,0.46,0.00}
\definecolor{goldenrod1}{rgb}{1.00,0.76,0.15}
\definecolor{goldenrod2}{rgb}{0.93,0.71,0.13}
\definecolor{goldenrod3}{rgb}{0.80,0.61,0.11}
\definecolor{goldenrod4}{rgb}{0.55,0.41,0.08}
\definecolor{goldenrod}{rgb}{0.85,0.65,0.13}
\definecolor{gold}{rgb}{1.00,0.84,0.00}
\definecolor{gray0}{rgb}{0.00,0.00,0.00}
\definecolor{gray100}{rgb}{1.00,1.00,1.00}
\definecolor{gray10}{rgb}{0.10,0.10,0.10}
\definecolor{gray11}{rgb}{0.11,0.11,0.11}
\definecolor{gray12}{rgb}{0.12,0.12,0.12}
\definecolor{gray13}{rgb}{0.13,0.13,0.13}
\definecolor{gray14}{rgb}{0.14,0.14,0.14}
\definecolor{gray15}{rgb}{0.15,0.15,0.15}
\definecolor{gray16}{rgb}{0.16,0.16,0.16}
\definecolor{gray17}{rgb}{0.17,0.17,0.17}
\definecolor{gray18}{rgb}{0.18,0.18,0.18}
\definecolor{gray19}{rgb}{0.19,0.19,0.19}
\definecolor{gray1}{rgb}{0.01,0.01,0.01}
\definecolor{gray20}{rgb}{0.20,0.20,0.20}
\definecolor{gray21}{rgb}{0.21,0.21,0.21}
\definecolor{gray22}{rgb}{0.22,0.22,0.22}
\definecolor{gray23}{rgb}{0.23,0.23,0.23}
\definecolor{gray24}{rgb}{0.24,0.24,0.24}
\definecolor{gray25}{rgb}{0.25,0.25,0.25}
\definecolor{gray26}{rgb}{0.26,0.26,0.26}
\definecolor{gray27}{rgb}{0.27,0.27,0.27}
\definecolor{gray28}{rgb}{0.28,0.28,0.28}
\definecolor{gray29}{rgb}{0.29,0.29,0.29}
\definecolor{gray2}{rgb}{0.02,0.02,0.02}
\definecolor{gray30}{rgb}{0.30,0.30,0.30}
\definecolor{gray31}{rgb}{0.31,0.31,0.31}
\definecolor{gray32}{rgb}{0.32,0.32,0.32}
\definecolor{gray33}{rgb}{0.33,0.33,0.33}
\definecolor{gray34}{rgb}{0.34,0.34,0.34}
\definecolor{gray35}{rgb}{0.35,0.35,0.35}
\definecolor{gray36}{rgb}{0.36,0.36,0.36}
\definecolor{gray37}{rgb}{0.37,0.37,0.37}
\definecolor{gray38}{rgb}{0.38,0.38,0.38}
\definecolor{gray39}{rgb}{0.39,0.39,0.39}
\definecolor{gray3}{rgb}{0.03,0.03,0.03}
\definecolor{gray40}{rgb}{0.40,0.40,0.40}
\definecolor{gray41}{rgb}{0.41,0.41,0.41}
\definecolor{gray42}{rgb}{0.42,0.42,0.42}
\definecolor{gray43}{rgb}{0.43,0.43,0.43}
\definecolor{gray44}{rgb}{0.44,0.44,0.44}
\definecolor{gray45}{rgb}{0.45,0.45,0.45}
\definecolor{gray46}{rgb}{0.46,0.46,0.46}
\definecolor{gray47}{rgb}{0.47,0.47,0.47}
\definecolor{gray48}{rgb}{0.48,0.48,0.48}
\definecolor{gray49}{rgb}{0.49,0.49,0.49}
\definecolor{gray4}{rgb}{0.04,0.04,0.04}
\definecolor{gray50}{rgb}{0.50,0.50,0.50}
\definecolor{gray51}{rgb}{0.51,0.51,0.51}
\definecolor{gray52}{rgb}{0.52,0.52,0.52}
\definecolor{gray53}{rgb}{0.53,0.53,0.53}
\definecolor{gray54}{rgb}{0.54,0.54,0.54}
\definecolor{gray55}{rgb}{0.55,0.55,0.55}
\definecolor{gray56}{rgb}{0.56,0.56,0.56}
\definecolor{gray57}{rgb}{0.57,0.57,0.57}
\definecolor{gray58}{rgb}{0.58,0.58,0.58}
\definecolor{gray59}{rgb}{0.59,0.59,0.59}
\definecolor{gray5}{rgb}{0.05,0.05,0.05}
\definecolor{gray60}{rgb}{0.60,0.60,0.60}
\definecolor{gray61}{rgb}{0.61,0.61,0.61}
\definecolor{gray62}{rgb}{0.62,0.62,0.62}
\definecolor{gray63}{rgb}{0.63,0.63,0.63}
\definecolor{gray64}{rgb}{0.64,0.64,0.64}
\definecolor{gray65}{rgb}{0.65,0.65,0.65}
\definecolor{gray66}{rgb}{0.66,0.66,0.66}
\definecolor{gray67}{rgb}{0.67,0.67,0.67}
\definecolor{gray68}{rgb}{0.68,0.68,0.68}
\definecolor{gray69}{rgb}{0.69,0.69,0.69}
\definecolor{gray6}{rgb}{0.06,0.06,0.06}
\definecolor{gray70}{rgb}{0.70,0.70,0.70}
\definecolor{gray71}{rgb}{0.71,0.71,0.71}
\definecolor{gray72}{rgb}{0.72,0.72,0.72}
\definecolor{gray73}{rgb}{0.73,0.73,0.73}
\definecolor{gray74}{rgb}{0.74,0.74,0.74}
\definecolor{gray75}{rgb}{0.75,0.75,0.75}
\definecolor{gray76}{rgb}{0.76,0.76,0.76}
\definecolor{gray77}{rgb}{0.77,0.77,0.77}
\definecolor{gray78}{rgb}{0.78,0.78,0.78}
\definecolor{gray79}{rgb}{0.79,0.79,0.79}
\definecolor{gray7}{rgb}{0.07,0.07,0.07}
\definecolor{gray80}{rgb}{0.80,0.80,0.80}
\definecolor{gray81}{rgb}{0.81,0.81,0.81}
\definecolor{gray82}{rgb}{0.82,0.82,0.82}
\definecolor{gray83}{rgb}{0.83,0.83,0.83}
\definecolor{gray84}{rgb}{0.84,0.84,0.84}
\definecolor{gray85}{rgb}{0.85,0.85,0.85}
\definecolor{gray86}{rgb}{0.86,0.86,0.86}
\definecolor{gray87}{rgb}{0.87,0.87,0.87}
\definecolor{gray88}{rgb}{0.88,0.88,0.88}
\definecolor{gray89}{rgb}{0.89,0.89,0.89}
\definecolor{gray8}{rgb}{0.08,0.08,0.08}
\definecolor{gray90}{rgb}{0.90,0.90,0.90}
\definecolor{gray91}{rgb}{0.91,0.91,0.91}
\definecolor{gray92}{rgb}{0.92,0.92,0.92}
\definecolor{gray93}{rgb}{0.93,0.93,0.93}
\definecolor{gray94}{rgb}{0.94,0.94,0.94}
\definecolor{gray95}{rgb}{0.95,0.95,0.95}
\definecolor{gray96}{rgb}{0.96,0.96,0.96}
\definecolor{gray97}{rgb}{0.97,0.97,0.97}
\definecolor{gray98}{rgb}{0.98,0.98,0.98}
\definecolor{gray99}{rgb}{0.99,0.99,0.99}
\definecolor{gray9}{rgb}{0.09,0.09,0.09}
\definecolor{gray}{rgb}{0.75,0.75,0.75}
\definecolor{green1}{rgb}{0.00,1.00,0.00}
\definecolor{green2}{rgb}{0.00,0.93,0.00}
\definecolor{green3}{rgb}{0.00,0.80,0.00}
\definecolor{green4}{rgb}{0.00,0.55,0.00}
\definecolor{greenyellow}{rgb}{0.68,1.00,0.18}
\definecolor{green}{rgb}{0.00,1.00,0.00}
\definecolor{grey0}{rgb}{0.00,0.00,0.00}
\definecolor{grey100}{rgb}{1.00,1.00,1.00}
\definecolor{grey10}{rgb}{0.10,0.10,0.10}
\definecolor{grey11}{rgb}{0.11,0.11,0.11}
\definecolor{grey12}{rgb}{0.12,0.12,0.12}
\definecolor{grey13}{rgb}{0.13,0.13,0.13}
\definecolor{grey14}{rgb}{0.14,0.14,0.14}
\definecolor{grey15}{rgb}{0.15,0.15,0.15}
\definecolor{grey16}{rgb}{0.16,0.16,0.16}
\definecolor{grey17}{rgb}{0.17,0.17,0.17}
\definecolor{grey18}{rgb}{0.18,0.18,0.18}
\definecolor{grey19}{rgb}{0.19,0.19,0.19}
\definecolor{grey1}{rgb}{0.01,0.01,0.01}
\definecolor{grey20}{rgb}{0.20,0.20,0.20}
\definecolor{grey21}{rgb}{0.21,0.21,0.21}
\definecolor{grey22}{rgb}{0.22,0.22,0.22}
\definecolor{grey23}{rgb}{0.23,0.23,0.23}
\definecolor{grey24}{rgb}{0.24,0.24,0.24}
\definecolor{grey25}{rgb}{0.25,0.25,0.25}
\definecolor{grey26}{rgb}{0.26,0.26,0.26}
\definecolor{grey27}{rgb}{0.27,0.27,0.27}
\definecolor{grey28}{rgb}{0.28,0.28,0.28}
\definecolor{grey29}{rgb}{0.29,0.29,0.29}
\definecolor{grey2}{rgb}{0.02,0.02,0.02}
\definecolor{grey30}{rgb}{0.30,0.30,0.30}
\definecolor{grey31}{rgb}{0.31,0.31,0.31}
\definecolor{grey32}{rgb}{0.32,0.32,0.32}
\definecolor{grey33}{rgb}{0.33,0.33,0.33}
\definecolor{grey34}{rgb}{0.34,0.34,0.34}
\definecolor{grey35}{rgb}{0.35,0.35,0.35}
\definecolor{grey36}{rgb}{0.36,0.36,0.36}
\definecolor{grey37}{rgb}{0.37,0.37,0.37}
\definecolor{grey38}{rgb}{0.38,0.38,0.38}
\definecolor{grey39}{rgb}{0.39,0.39,0.39}
\definecolor{grey3}{rgb}{0.03,0.03,0.03}
\definecolor{grey40}{rgb}{0.40,0.40,0.40}
\definecolor{grey41}{rgb}{0.41,0.41,0.41}
\definecolor{grey42}{rgb}{0.42,0.42,0.42}
\definecolor{grey43}{rgb}{0.43,0.43,0.43}
\definecolor{grey44}{rgb}{0.44,0.44,0.44}
\definecolor{grey45}{rgb}{0.45,0.45,0.45}
\definecolor{grey46}{rgb}{0.46,0.46,0.46}
\definecolor{grey47}{rgb}{0.47,0.47,0.47}
\definecolor{grey48}{rgb}{0.48,0.48,0.48}
\definecolor{grey49}{rgb}{0.49,0.49,0.49}
\definecolor{grey4}{rgb}{0.04,0.04,0.04}
\definecolor{grey50}{rgb}{0.50,0.50,0.50}
\definecolor{grey51}{rgb}{0.51,0.51,0.51}
\definecolor{grey52}{rgb}{0.52,0.52,0.52}
\definecolor{grey53}{rgb}{0.53,0.53,0.53}
\definecolor{grey54}{rgb}{0.54,0.54,0.54}
\definecolor{grey55}{rgb}{0.55,0.55,0.55}
\definecolor{grey56}{rgb}{0.56,0.56,0.56}
\definecolor{grey57}{rgb}{0.57,0.57,0.57}
\definecolor{grey58}{rgb}{0.58,0.58,0.58}
\definecolor{grey59}{rgb}{0.59,0.59,0.59}
\definecolor{grey5}{rgb}{0.05,0.05,0.05}
\definecolor{grey60}{rgb}{0.60,0.60,0.60}
\definecolor{grey61}{rgb}{0.61,0.61,0.61}
\definecolor{grey62}{rgb}{0.62,0.62,0.62}
\definecolor{grey63}{rgb}{0.63,0.63,0.63}
\definecolor{grey64}{rgb}{0.64,0.64,0.64}
\definecolor{grey65}{rgb}{0.65,0.65,0.65}
\definecolor{grey66}{rgb}{0.66,0.66,0.66}
\definecolor{grey67}{rgb}{0.67,0.67,0.67}
\definecolor{grey68}{rgb}{0.68,0.68,0.68}
\definecolor{grey69}{rgb}{0.69,0.69,0.69}
\definecolor{grey6}{rgb}{0.06,0.06,0.06}
\definecolor{grey70}{rgb}{0.70,0.70,0.70}
\definecolor{grey71}{rgb}{0.71,0.71,0.71}
\definecolor{grey72}{rgb}{0.72,0.72,0.72}
\definecolor{grey73}{rgb}{0.73,0.73,0.73}
\definecolor{grey74}{rgb}{0.74,0.74,0.74}
\definecolor{grey75}{rgb}{0.75,0.75,0.75}
\definecolor{grey76}{rgb}{0.76,0.76,0.76}
\definecolor{grey77}{rgb}{0.77,0.77,0.77}
\definecolor{grey78}{rgb}{0.78,0.78,0.78}
\definecolor{grey79}{rgb}{0.79,0.79,0.79}
\definecolor{grey7}{rgb}{0.07,0.07,0.07}
\definecolor{grey80}{rgb}{0.80,0.80,0.80}
\definecolor{grey81}{rgb}{0.81,0.81,0.81}
\definecolor{grey82}{rgb}{0.82,0.82,0.82}
\definecolor{grey83}{rgb}{0.83,0.83,0.83}
\definecolor{grey84}{rgb}{0.84,0.84,0.84}
\definecolor{grey85}{rgb}{0.85,0.85,0.85}
\definecolor{grey86}{rgb}{0.86,0.86,0.86}
\definecolor{grey87}{rgb}{0.87,0.87,0.87}
\definecolor{grey88}{rgb}{0.88,0.88,0.88}
\definecolor{grey89}{rgb}{0.89,0.89,0.89}
\definecolor{grey8}{rgb}{0.08,0.08,0.08}
\definecolor{grey90}{rgb}{0.90,0.90,0.90}
\definecolor{grey91}{rgb}{0.91,0.91,0.91}
\definecolor{grey92}{rgb}{0.92,0.92,0.92}
\definecolor{grey93}{rgb}{0.93,0.93,0.93}
\definecolor{grey94}{rgb}{0.94,0.94,0.94}
\definecolor{grey95}{rgb}{0.95,0.95,0.95}
\definecolor{grey96}{rgb}{0.96,0.96,0.96}
\definecolor{grey97}{rgb}{0.97,0.97,0.97}
\definecolor{grey98}{rgb}{0.98,0.98,0.98}
\definecolor{grey99}{rgb}{0.99,0.99,0.99}
\definecolor{grey9}{rgb}{0.09,0.09,0.09}
\definecolor{grey}{rgb}{0.75,0.75,0.75}
\definecolor{honeydew1}{rgb}{0.94,1.00,0.94}
\definecolor{honeydew2}{rgb}{0.88,0.93,0.88}
\definecolor{honeydew3}{rgb}{0.76,0.80,0.76}
\definecolor{honeydew4}{rgb}{0.51,0.55,0.51}
\definecolor{honeydew}{rgb}{0.94,1.00,0.94}
\definecolor{hotpink}{rgb}{1.00,0.41,0.71}
\definecolor{indianred}{rgb}{0.80,0.36,0.36}
\definecolor{ivory1}{rgb}{1.00,1.00,0.94}
\definecolor{ivory2}{rgb}{0.93,0.93,0.88}
\definecolor{ivory3}{rgb}{0.80,0.80,0.76}
\definecolor{ivory4}{rgb}{0.55,0.55,0.51}
\definecolor{ivory}{rgb}{1.00,1.00,0.94}
\definecolor{khaki1}{rgb}{1.00,0.96,0.56}
\definecolor{khaki2}{rgb}{0.93,0.90,0.52}
\definecolor{khaki3}{rgb}{0.80,0.78,0.45}
\definecolor{khaki4}{rgb}{0.55,0.53,0.31}
\definecolor{khaki}{rgb}{0.94,0.90,0.55}
\definecolor{lavenderblush}{rgb}{1.00,0.94,0.96}
\definecolor{lavender}{rgb}{0.90,0.90,0.98}
\definecolor{lawngreen}{rgb}{0.49,0.99,0.00}
\definecolor{lemonchiffon}{rgb}{1.00,0.98,0.80}
\definecolor{lightblue}{rgb}{0.68,0.85,0.90}
\definecolor{lightcoral}{rgb}{0.94,0.50,0.50}
\definecolor{lightcyan}{rgb}{0.88,1.00,1.00}
\definecolor{lightgoldenrod}{rgb}{0.93,0.87,0.51}
\definecolor{lightgoldenrod}{rgb}{0.98,0.98,0.82}
\definecolor{lightgray}{rgb}{0.83,0.83,0.83}
\definecolor{lightgreen}{rgb}{0.56,0.93,0.56}
\definecolor{lightgrey}{rgb}{0.83,0.83,0.83}
\definecolor{lightpink}{rgb}{1.00,0.71,0.76}
\definecolor{lightsalmon}{rgb}{1.00,0.63,0.48}
\definecolor{lightsea}{rgb}{0.13,0.70,0.67}
\definecolor{lightsky}{rgb}{0.53,0.81,0.98}
\definecolor{lightslate}{rgb}{0.47,0.53,0.60}
\definecolor{lightslate}{rgb}{0.47,0.53,0.60}
\definecolor{lightslate}{rgb}{0.52,0.44,1.00}
\definecolor{lightsteel}{rgb}{0.69,0.77,0.87}
\definecolor{lightyellow}{rgb}{1.00,1.00,0.88}
\definecolor{limegreen}{rgb}{0.20,0.80,0.20}
\definecolor{linen}{rgb}{0.98,0.94,0.90}
\definecolor{magenta1}{rgb}{1.00,0.00,1.00}
\definecolor{magenta2}{rgb}{0.93,0.00,0.93}
\definecolor{magenta3}{rgb}{0.80,0.00,0.80}
\definecolor{magenta4}{rgb}{0.55,0.00,0.55}
\definecolor{magenta}{rgb}{1.00,0.00,1.00}
\definecolor{maroon1}{rgb}{1.00,0.20,0.70}
\definecolor{maroon2}{rgb}{0.93,0.19,0.65}
\definecolor{maroon3}{rgb}{0.80,0.16,0.56}
\definecolor{maroon4}{rgb}{0.55,0.11,0.38}
\definecolor{maroon}{rgb}{0.69,0.19,0.38}
\definecolor{mediumaquamarine}{rgb}{0.40,0.80,0.67}
\definecolor{mediumblue}{rgb}{0.00,0.00,0.80}
\definecolor{mediumorchid}{rgb}{0.73,0.33,0.83}
\definecolor{mediumpurple}{rgb}{0.58,0.44,0.86}
\definecolor{mediumsea}{rgb}{0.24,0.70,0.44}
\definecolor{mediumslate}{rgb}{0.48,0.41,0.93}
\definecolor{mediumspring}{rgb}{0.00,0.98,0.60}
\definecolor{mediumturquoise}{rgb}{0.28,0.82,0.80}
\definecolor{mediumviolet}{rgb}{0.78,0.08,0.52}
\definecolor{midnightblue}{rgb}{0.10,0.10,0.44}
\definecolor{mintcream}{rgb}{0.96,1.00,0.98}
\definecolor{mistyrose}{rgb}{1.00,0.89,0.88}
\definecolor{moccasin}{rgb}{1.00,0.89,0.71}
\definecolor{navajowhite}{rgb}{1.00,0.87,0.68}
\definecolor{navyblue}{rgb}{0.00,0.00,0.50}
\definecolor{navy}{rgb}{0.00,0.00,0.50}
\definecolor{oldlace}{rgb}{0.99,0.96,0.90}
\definecolor{olivedrab}{rgb}{0.42,0.56,0.14}
\definecolor{orange1}{rgb}{1.00,0.65,0.00}
\definecolor{orange2}{rgb}{0.93,0.60,0.00}
\definecolor{orange3}{rgb}{0.80,0.52,0.00}
\definecolor{orange4}{rgb}{0.55,0.35,0.00}
\definecolor{orangered}{rgb}{1.00,0.27,0.00}
\definecolor{orange}{rgb}{1.00,0.65,0.00}
\definecolor{orchid1}{rgb}{1.00,0.51,0.98}
\definecolor{orchid2}{rgb}{0.93,0.48,0.91}
\definecolor{orchid3}{rgb}{0.80,0.41,0.79}
\definecolor{orchid4}{rgb}{0.55,0.28,0.54}
\definecolor{orchid}{rgb}{0.85,0.44,0.84}
\definecolor{palegoldenrod}{rgb}{0.93,0.91,0.67}
\definecolor{palegreen}{rgb}{0.60,0.98,0.60}
\definecolor{paleturquoise}{rgb}{0.69,0.93,0.93}
\definecolor{paleviolet}{rgb}{0.86,0.44,0.58}
\definecolor{papayawhip}{rgb}{1.00,0.94,0.84}
\definecolor{peachpuff}{rgb}{1.00,0.85,0.73}
\definecolor{peru}{rgb}{0.80,0.52,0.25}
\definecolor{pink1}{rgb}{1.00,0.71,0.77}
\definecolor{pink2}{rgb}{0.93,0.66,0.72}
\definecolor{pink3}{rgb}{0.80,0.57,0.62}
\definecolor{pink4}{rgb}{0.55,0.39,0.42}
\definecolor{pink}{rgb}{1.00,0.75,0.80}
\definecolor{plum1}{rgb}{1.00,0.73,1.00}
\definecolor{plum2}{rgb}{0.93,0.68,0.93}
\definecolor{plum3}{rgb}{0.80,0.59,0.80}
\definecolor{plum4}{rgb}{0.55,0.40,0.55}
\definecolor{plum}{rgb}{0.87,0.63,0.87}
\definecolor{powderblue}{rgb}{0.69,0.88,0.90}
\definecolor{purple1}{rgb}{0.61,0.19,1.00}
\definecolor{purple2}{rgb}{0.57,0.17,0.93}
\definecolor{purple3}{rgb}{0.49,0.15,0.80}
\definecolor{purple4}{rgb}{0.33,0.10,0.55}
\definecolor{purple}{rgb}{0.63,0.13,0.94}
\definecolor{red1}{rgb}{1.00,0.00,0.00}
\definecolor{red2}{rgb}{0.93,0.00,0.00}
\definecolor{red3}{rgb}{0.80,0.00,0.00}
\definecolor{red4}{rgb}{0.55,0.00,0.00}
\definecolor{red}{rgb}{1.00,0.00,0.00}
\definecolor{rosybrown}{rgb}{0.74,0.56,0.56}
\definecolor{royalblue}{rgb}{0.25,0.41,0.88}
\definecolor{saddlebrown}{rgb}{0.55,0.27,0.07}
\definecolor{salmon1}{rgb}{1.00,0.55,0.41}
\definecolor{salmon2}{rgb}{0.93,0.51,0.38}
\definecolor{salmon3}{rgb}{0.80,0.44,0.33}
\definecolor{salmon4}{rgb}{0.55,0.30,0.22}
\definecolor{salmon}{rgb}{0.98,0.50,0.45}
\definecolor{sandybrown}{rgb}{0.96,0.64,0.38}
\definecolor{seagreen}{rgb}{0.18,0.55,0.34}
\definecolor{seashell1}{rgb}{1.00,0.96,0.93}
\definecolor{seashell2}{rgb}{0.93,0.90,0.87}
\definecolor{seashell3}{rgb}{0.80,0.77,0.75}
\definecolor{seashell4}{rgb}{0.55,0.53,0.51}
\definecolor{seashell}{rgb}{1.00,0.96,0.93}
\definecolor{sienna1}{rgb}{1.00,0.51,0.28}
\definecolor{sienna2}{rgb}{0.93,0.47,0.26}
\definecolor{sienna3}{rgb}{0.80,0.41,0.22}
\definecolor{sienna4}{rgb}{0.55,0.28,0.15}
\definecolor{sienna}{rgb}{0.63,0.32,0.18}
\definecolor{skyblue}{rgb}{0.53,0.81,0.92}
\definecolor{slateblue}{rgb}{0.42,0.35,0.80}
\definecolor{slategray}{rgb}{0.44,0.50,0.56}
\definecolor{slategrey}{rgb}{0.44,0.50,0.56}
\definecolor{snow1}{rgb}{1.00,0.98,0.98}
\definecolor{snow2}{rgb}{0.93,0.91,0.91}
\definecolor{snow3}{rgb}{0.80,0.79,0.79}
\definecolor{snow4}{rgb}{0.55,0.54,0.54}
\definecolor{snow}{rgb}{1.00,0.98,0.98}
\definecolor{springgreen}{rgb}{0.00,1.00,0.50}
\definecolor{steelblue}{rgb}{0.27,0.51,0.71}
\definecolor{tan1}{rgb}{1.00,0.65,0.31}
\definecolor{tan2}{rgb}{0.93,0.60,0.29}
\definecolor{tan3}{rgb}{0.80,0.52,0.25}
\definecolor{tan4}{rgb}{0.55,0.35,0.17}
\definecolor{tan}{rgb}{0.82,0.71,0.55}
\definecolor{thistle1}{rgb}{1.00,0.88,1.00}
\definecolor{thistle2}{rgb}{0.93,0.82,0.93}
\definecolor{thistle3}{rgb}{0.80,0.71,0.80}
\definecolor{thistle4}{rgb}{0.55,0.48,0.55}
\definecolor{thistle}{rgb}{0.85,0.75,0.85}
\definecolor{tomato1}{rgb}{1.00,0.39,0.28}
\definecolor{tomato2}{rgb}{0.93,0.36,0.26}
\definecolor{tomato3}{rgb}{0.80,0.31,0.22}
\definecolor{tomato4}{rgb}{0.55,0.21,0.15}
\definecolor{tomato}{rgb}{1.00,0.39,0.28}
\definecolor{turquoise1}{rgb}{0.00,0.96,1.00}
\definecolor{turquoise2}{rgb}{0.00,0.90,0.93}
\definecolor{turquoise3}{rgb}{0.00,0.77,0.80}
\definecolor{turquoise4}{rgb}{0.00,0.53,0.55}
\definecolor{turquoise}{rgb}{0.25,0.88,0.82}
\definecolor{violetred}{rgb}{0.82,0.13,0.56}
\definecolor{violet}{rgb}{0.93,0.51,0.93}
\definecolor{wheat1}{rgb}{1.00,0.91,0.73}
\definecolor{wheat2}{rgb}{0.93,0.85,0.68}
\definecolor{wheat3}{rgb}{0.80,0.73,0.59}
\definecolor{wheat4}{rgb}{0.55,0.49,0.40}
\definecolor{wheat}{rgb}{0.96,0.87,0.70}
\definecolor{whitesmoke}{rgb}{0.96,0.96,0.96}
\definecolor{white}{rgb}{1.00,1.00,1.00}
\definecolor{yellow1}{rgb}{1.00,1.00,0.00}
\definecolor{yellow2}{rgb}{0.93,0.93,0.00}
\definecolor{yellow3}{rgb}{0.80,0.80,0.00}
\definecolor{yellow4}{rgb}{0.55,0.55,0.00}
\definecolor{yellowgreen}{rgb}{0.60,0.80,0.20}
\definecolor{yellow}{rgb}{1.00,1.00,0.00}
\usepackage{caption}
\newtheorem{definition}{Definition}[section]

\newtheorem{theorem}[definition]{Theorem}
\newtheorem{proposition}[definition]{Proposition}

\theoremstyle{remark}
\newtheorem{remark}[definition]{Remark}

\newcommand{\ie}{; {\it i.e., }}
\newcommand{\e}{\varepsilon}

\font\tenmsb=msbm10
\font\sevenmsb=msbm7
\font\fivemsb=msbm5
\newfam\msbfam
\textfont\msbfam=\tenmsb
\scriptfont\msbfam=\sevenmsb
\scriptscriptfont\msbfam=\fivemsb
\def\Bbb#1{{\fam\msbfam\relax#1}}

\def\R{\Bbb R}

\newcommand{\NN}{{\Bbb N}}

\newcommand{\ZZ}{{\Bbb Z}}

\def\to{\rightarrow}

\def\xie{{x^\e_i}}
\def\Txie{{ x_i^{\e,T}}}

\begin{document}
\title{Minimising movements for oscillating energies: \\
the critical regime}
\author[1]{Nadia Ansini}
\author[2]{Andrea Braides}
\author[3]{Johannes Zimmer}
\affil[1]{Dept. of Mathematics,
 Sapienza University of Rome, P.le Aldo Moro 2, 00185 Rome, Italy}
\affil[2]{Dept. of Mathematics,
University of Rome `Tor Vergata',  Via della Ricerca Scientifica, 00133 Rome, Italy}
\affil[3]{Dept. of Mathematical Sciences, University of Bath, Bath BA2 7AY, UK} 

\date{}

\maketitle

\begin{abstract}
  Minimising movements are investigated for an energy which is the superposition of a convex functional and fast small
  oscillations. Thus a minimising movement scheme involves a temporal parameter $\tau$ and a spatial parameter $\e$,
  with $\tau$ describing the time step and the frequency of the oscillations being proportional to $\frac 1 \e$. The
  extreme cases of fast time scales $\tau\ll\e$ and slow time scales $\e\ll\tau$ have been investigated
  in~\cite{Braides2014a}. In this article, the intermediate (critical) case of finite ratio $\e/\tau>0$ is studied. It
  is shown that a pinning threshold exists, with initial data below the threshold being a fixed point of the
  dynamics. A characterisation of the pinning threshold is given. For initial data above the pinning threshold, the
  equation and velocity describing the homogenised motion are determined.


\end{abstract}

\section{Introduction}
\label{sec:Introduction}

In this paper we analyse a minimising-movement approach for gradient flows with wiggly energies,
\begin{equation}
  \label{eq:gf}
    x'(t) = - \mu \frac{\partial E_\e(x(t))}{\partial x}. 
\end{equation}
A prototypical model of the energy is an oscillating perturbation of a quadratic energy,
\begin{equation}
  \label{eq:quaden}
  E_\e (x)= {1\over 2}x^2 - \e \cos \Bigl({x\over \e}\Bigr).
\end{equation}
This mathematical problem can be motivated by the analysis of interface motion in materials science. There is a range
of problems where interfaces form in a specimen and propagate so that a material particle crossing the interface
changes its stability, by transforming from an unstable or metastable state to a more stable one; see,
e.g.,~\cite{Abeyaratne1999a}. Often this evolution is load-driven, in the sense that an applied load enables a particle
to explore states of lower energy. Let us consider an interface, say between twin boundaries or phase boundaries,
macroscopically propagating with some velocity $v$. However, microscopically the interface typically does not move
homogeneously. Instead, the interface tends to propagate forward as a whole by a series of incremental steps. To
illustrate this, let us picture an interface consisting of a straight horizontal line segment, then a step up, followed
by another horizontal line segment, moving up towards a more stable state. Then it is normally advantageous for the
interface to propagate the step sideways, i.e., move upward the particle next to the step, and then move the remaining
particles consecutively. This leads to a multi-welled energy landscape, with local minima spaced periodically with high
frequency. One model for the propagation of an interface in this manner is the Frenkel-Kontorova chain with forcing,
where the motion of atom $n$ is~\cite{Abeyaratne1999a}
\begin{equation*}
  m u_n''(t) = k (u_{n+1}(t) - 2 u_n(t) + u_{n-1}(t)) - W'(u_n(t)) + f(n,t), 
\end{equation*}
or the continuum version
\begin{equation*}
  m  u''(x,t) = k u_{xx}(x,t) - W'(u(x,t)) + f(x,t). 
\end{equation*}
The model considered here can be interpreted as an unforced ($f=0$) case, where the kinetics is replaced by a simpler
(gradient flow) dynamics.

We remark that the same equation appears in a related but different context in Materials Science, again originating
from transition layers. Martensitic materials can form needles of phases with pronounced tips (see, e.g., photographs
in~\cite{Abeyaratne1996a}). During creep tests, it is observed that the volume fraction of the phase fractions involved
changes rather abruptly, and it is shown that this sudden change can be attributed to a sudden split of a tip into two
tips~\cite{Abeyaratne1996a}. One can picture this as a lenticular domain of one variant trying to grow; this growth
then occurs where the tip of the lens meets a boundary between twins, and fattening of the phase happens via splitting
of the tip in two and more tips. The splitting of a needle can then be attributed to a metastable transition, moving
from one local minimum to another one. This suggests a small-scale landscape with many minima, and the energy studied
by Abeyaratne, Chu and James~\cite{Abeyaratne1996a} is a macroscopic energy augmented by small-scale oscillations
$a \e \cos(\frac x \e)$, as studied here. In addition, the kinetic law in~\cite{Abeyaratne1996a} is taken to be a
gradient flow. Specifically, there it is shown that the solution $x_\e$ to the evolution equation~\eqref{eq:gf}
converges uniformly in time to the solution of
\begin{equation*}
   x(t) = - \mu \frac{\partial \bar E(x(t))}{\partial x}, \text{ with } x(0) = x_0, 
\end{equation*}
with an explicitly computed driving force $\frac{\partial \bar E}{\partial x}$. This latter system is then investigated
numerically. The variational analysis carried out here can be interpreted in this light. We consider time
discretisations, as numerical algorithms would employ, but on the level of the original (not homogenised) energy
$E_\e$, rather than $\bar E$. This leads to two parameters, the time discretisation $\tau$ and oscillation scale
$\e$. The different scaling regimes that follow naturally are analysed in this paper.

We make the trivial but important remark that the limit of a sequence of gradient flows associated to an family $E_\e$
is in general not the gradient flow of the limit of the energy. For example, on one hand for the
energy~\eqref{eq:quaden}, the associated gradient flow is ($\mu=1$)
\begin{equation*}
  x_\e'(t)=-x_\e(t)-\sin \Bigl({x_\e(t)\over \e}\Bigr),
\end{equation*}
for initial datum $x_0$. If $x_0\in (-1,1)$, then such solutions are trapped between stationary solutions, and they
converge to the trivial constant state $x_0$ (\emph{pinning}), while if $|x_0|\ge 1$ they can be shown to converge to a
solution $x$ of the gradient flow
\begin{equation*}
  x'(t)=- {\rm sign}\,{x(t)}\sqrt{x^2(t)-1}.
\end{equation*}
On the other hand, $E_\e$ converge uniformly to the quadratic energy, whose gradient flow is trivially
\begin{equation*}
  x'(t)=-x(t).
\end{equation*}

These behaviours can be obtained as limit cases of \emph{minimising movements along the sequence of energies $E_\e$} at
different \emph{time scales}. Minimising movements are defined as follows: with fixed $\e$ (the spatial scale) and
$\tau$ (the time scale), we set $x^{\e,\tau}_0= x_0$ and choose recursively $x^{\e,\tau}_k$ as a minimiser of
\begin{equation*}
  x\mapsto E_\e(x)+{1\over 2\tau}|x-x^{\e,\tau}_{k-1}|^2.
\end{equation*}
This process gives the piecewise-constant trajectories 
\begin{equation*}
  x^{\e,\tau}(t)=x^{\e,\tau}_{\lfloor t/\tau\rfloor}.
\end{equation*}

With fixed $\tau=\tau(\e)$, a minimising movement $x(t)$ along the sequence of energies $E_\e$ at time scale $\tau$ is defined as any limit of subsequences of $x^{\e,\tau}(t)$. Simple examples show that the limit may indeed depend on the subsequence and on the choice of $\tau$. If $E_\e$ is independent of $\e$ this notion coincides with the one given by De Giorgi~\cite{De-Giorgi1993a} and at the basis of modern notions of gradient flows (see the monograph by Ambrosio, Gigli and Savar\'e~\cite{Ambrosio2008b}). Examples of problems related to varying $E_\e$ are analysed in~\cite{Braides2014a}, where in particular it is shown that for the energies above and for $\tau\ll\e$ (\emph{fast time scale}), the minimising movement $x$ coincides with the limit of the solutions $x_\e$ of the gradient flows at fixed $\e$, while for $\e\ll\tau$ (\emph{slow time scale}) it coincides with the gradient flow of the limit quadratic energy.  That observation highlights the existence of a \emph{critical time scaling} when $\tau\sim\e$, for which the minimising movements are not trivially described by the limit of gradient flows or the gradient flow of the limit. The behaviour at those scales is the object of this paper. A rather different very interesting line of investigation has been taken by Menon~\cite{Menon2002a}, and independently in parallel by Smyshlyaev. In~\cite{Menon2002a}, averaging techniques are developed in the context of the time-continuous dynamical system~\eqref{eq:gf}. The homogenisation of first-order ordinary differential equations, including error estimates, is studied further in~\cite{Ibrahim2010a}. We remark that the model we consider is deterministic, where the two parameters come from spatial oscillations and a time discretisation. For stochastic models, it is also natural to consider the effective behaviour in different scaling regimes of space and noise; we refer the reader to~\cite{Dupuis2012a}.

\paragraph{Plan of the paper} A summary of the results of the paper is as follows. We consider functions
$E_\e\colon \R\to \R$ given by
\begin{equation}
  \label{funzioni}
  E_\e (x)= h(x) + \e\, W \Bigl({x\over \e}\Bigr),
\end{equation}
where $h$ is a strictly convex function bounded from below and $W$ is a one-periodic even Lipschitz function. We
consider a time scale $\tau=\tau(\e)$ such that $\e/\tau$ converges to $\gamma>0$. Therefore, the present analysis
complements the recent one of~\cite{Braides2014a} {(see also~\cite{Braides2008a})}, where the cases
$\gamma \in \{0, +\infty\}$ are investigated. In terms of the mechanical problem of interface propagation discussed
above, we show that pinning will occur for small initial data, while large data leads to a gradient flow evolution for
which the averaged velocity can be computed. More precisely, we prove that in that case the unique minimising movement
$x^\gamma$ with initial datum $x_0>0$ (the case $x_0<0$ is analogous by symmetry) is described as follows:
\begin{enumerate}
\item \emph{Pinning threshold:\/} There exists $T_\gamma$ such that $x^\gamma(t)=x_0$ for all $t$ if
  $|x_0|\le T_\gamma$. The pinning threshold is characterised in Proposition~\ref{carac}.
\item \emph{Homogenised equation:\/} If $x_0>T_\gamma$ then $x^\gamma(t)$ is characterised as the non-increasing
  function satisfying
  \begin{equation*}
    {\rm d\over dt}x^\gamma(t)= -\gamma\, f_\gamma(h'(x^\gamma(t)))
  \end{equation*}
  at almost all $t>0$. The \emph{homogenised velocity} $f_\gamma(z)$ is the \emph{average velocity} (suitably defined)
  of any discrete orbit $\{y_k\}$ defined recursively by minimisation of the linearity problem
  \begin{equation*}
    y\mapsto z y+W(y)+{\gamma\over 2}(y-y_{k-1})^2,
  \end{equation*}
  which can be shown not to depend on the initial condition $y_0$.
\end{enumerate}
Mathematically, our analysis is confined to one space dimension, as it strongly relies on monotonicity properties
developed in Section~\ref{sec:Monot-behav-minim}.  A central argument is a comparison of solution to the nonlinear
energy as in~\eqref{funzioni}, and one where $h$, the non-oscillating part, is suitably linearised. This argument is
developed in Section~\ref{sec:Linearized-energy}.


\section{Minimising movement along a sequence}
\label{sec:Minim-movem-along}

We recall the general definition of minimising movements for a sequence of functionals defined on a Hilbert space.
\begin{definition}
  \label{DefMMseq}
  Let $X$ be a separable Hilbert space, $E_\e \colon X\to [0, +\infty)$ equicoercive and lower semicontinuous and
  $x_0^{\e}\to x_0$ with $E_\e (x_0^{\e}) \le C<+\infty$ and $\tau_\e >0$ converging to $0$ as $\e\to 0$. For fixed
  $\e>0$, we define recursively $x_i^{\e}$ as a minimiser of the problem
  \begin{equation}
    \label{minpMM}
    \min\Bigl\{ E_\e (x) + {1\over \tau_{\e}}  \Vert x- x_{i-1}^{\e}\Vert^2\Bigr\}\,,
  \end{equation}
  and the piecewise-constant trajectory 
  \begin{equation}
    \label{seqMM}
    x^\e (t):= x_i^\e= x^\e_{\lfloor t/\tau_{\e}\rfloor}\,, \quad t\in [i\tau_\e, (i+1)\tau_\e)\,.
  \end{equation}
  A \emph{minimising movement for $E_\e$ at time scale $\tau$} from $x_0^\e$ is the limit of a subsequence $x^{\e_j}$,
  \begin{equation*}
    x(t)= \lim_{j\to +\infty} x^{\e_j}(t)\,,
  \end{equation*}
  with respect to the uniform convergence on compact sets of $[0, +\infty)$.
\end{definition}
This definition is justified by the following compactness result~\cite[Proposition 7.1]{Braides2014a}.

\begin{proposition}
  For every $E_\e$ and $x_0^{\e}$ as above, there exist minimising movements for $E_\e$, from $x_0^{\e}$, with
  $x(t)\in C^{1/2} ([0, +\infty); X)$.
\end{proposition}
For a comprehensive study of minimising movements for a fixed $E=E_\e$ we refer to~\cite{Ambrosio2008b}, while a
detailed analysis of some of its applications can be found in~\cite{Braides2014a}.

\section{Monotone behaviour of minimising movements}
\label{sec:Monot-behav-minim}

In the sequel we will study minimising movements for the functions $E_\e\colon \R\to \R$ given by~\eqref{funzioni},
where $h$ is a strictly convex function bounded from below. It is not restrictive to suppose that $h\ge 0$, and that
$h$ attains its global minimum in $x=0$. Furthermore, we assume that $W$ is a one-periodic even Lipschitz function with
$\Vert W^\prime\Vert_{\infty}=1$, and that the average of $W$ is $0$. The two latter assumptions serve as normalisation
only and are not restrictive.

We observe the following simple monotonicity property.

 \begin{proposition}\label{monof}
   Given any functions $\phi\colon\R\to \R$ and $\psi\colon\R\to \R$, and $\beta>0$, for any $x,x'\in\R$, let
   $y,y'\in\R$ be minimisers of
  \begin{equation*}
    t\mapsto \phi(t)+ \beta (t- x)^2\,,\qquad
    t\mapsto \psi(t)+ \beta (t- x')^2,
  \end{equation*}
  respectively. Then
  \begin{equation}
    \label{phi-psi}
    \phi(y)-\phi(y') + \psi(y')-\psi(y)\le 2\beta (x-x')(y-y')\,.
  \end{equation}
  In particular, if $\psi=\phi$ and $x\le x'$ then $y\le y'$.
\end{proposition}

\proof 
By assumption
\begin{equation*}
  \phi(y)+ \beta (y- x)^2 \le\phi(y') + \beta (y'- x)^2,
  \quad
  \psi(y')+ \beta (y'- x')^2 \le\psi(y)+ \beta (y- x')^2.
\end{equation*}
Summing up the two inequalities and simplifying the terms on both sides we obtain~\ref{phi-psi}.  Moreover, if
$\psi=\phi$, then
\begin{equation*}
  (x-x')(y-y')\ge 0,
\end{equation*}
which yields the desired inequality. 
\qed

Before analysing the case of fixed ratio $\e/\tau$, we make some general remarks. According to
Definition~\ref{DefMMseq}, we define iteratively the global minimiser $x_{i+1}^\e$ to
\begin{equation}
  \label{F(i)general}
  F_{\e}(x, \xie)= h(x)  + \e \,W \Bigl({x\over \e}\Bigr) + {1\over 2{\tau_\e}} (x- \xie)^2\,.
\end{equation}
First of all, we observe that the sequence of minimisers $(x_{i}^\e)_i$ is monotone.

\begin{proposition}[Monotone behavior of $\xie$]
  \label{Mono-xie}
  Let $x_{i+1}^\e$ be a minimiser to~\eqref{F(i)general}. Then the following holds.
  \begin{enumerate}
  \item \label{it:monbeh1} If $x^\e_{i+1}\le \xie$, then $     x_{i+2}^\e\le x_{i+1}^\e$\,.
  \item \label{it:monbeh2} If $x^\e_{i+1}\ge \xie$, then $     x_{i+2}^\e\ge x_{i+1}^\e$\,.
  \end{enumerate}
  In particular, $t\mapsto x^\e(t)$ and $t\mapsto x(t)$ are monotone functions.
\end{proposition}

\proof \ref{it:monbeh1} and~\ref{it:monbeh2} are straightforward consequences of Proposition~\ref{monof} with
$\phi(t)=\psi(t)= h(t) +\e W(t/\e)$ and $\beta = 1/(2\tau_{\e})$.


By Definition~\ref{DefMMseq}, $x^\e(t)$ is monotone in $t$, and since it converges uniformly to
$x(t)\in C^{1/2}([0, +\infty))$, on compact sets of $[0, +\infty)$, we may conclude that $x(t)$ is also a monotone
function.  
\qed

\section{Linearised energy}
\label{sec:Linearized-energy}

In order to characterise the velocity $x^\prime(t_0)$ of the minimising movement scheme~\eqref{minpMM}, we study the
average velocity given by
\begin{equation}
  \label{eq:vel}
  {x_i^\e- x_0^\e\over i\tau_{\e}}\,.
\end{equation}
We assume, without loss of generality, that $t_0=0$ and $x(0)=x_0^\e$.  We consider a (partial) linearisation of the
problem given by
\begin{equation}
  \label{F(i)linear}
  F^T_{\e}(x, \Txie)= Tx  + \e\, W \Bigl({x\over \e}\Bigr) + {1\over 2\tau_\e} (x- \Txie)^2\,,
\end{equation}
where $x_0^{\e, T}:= x_0^\e$. The term $Tx$ represents the ``linear approximation'' of the potential $h$ around the
point $h(x^\e_0)$ up to translation by a constant that does not depend on $T$ and $i$.  We recall that $h$ is a
strictly convex function, hence $h^\prime$ is a monotone increasing function.

\begin{proposition}[Monotone behavior of $\Txie$ with respect to $\xie$] 
  \label{monotonicity}
  Let $\delta>0$ be such that $h^\prime(x_0\pm\delta)$ exists.
  \begin{enumerate}
  \item Let $x^{\e,T}_{i}$ be the minimiser to~\eqref{F(i)linear}. Then \label{it:mon1}
    \begin{itemize}
    \item if $x^{\e,T}_{i+1}\le \Txie$, then  $x_{i+2}^{\e,T}\le x_{i+1}^{\e,T}$
    \item if $x^{\e, T}_{i+1}\ge \Txie$, then $x_{i+2}^{\e,T}\ge x_{i+1}^{\e,T}$\,.
    \end{itemize}
  \item Let $T= T^{\delta^+}:= h^\prime(x_0+\delta)$, and $x^\e_i$ minimiser of~\eqref{F(i)general}. Then if
    $x^{\e,T^{\delta^+}}_{i}\le x_{i}^\e$, then $ x^{\e,T^{\delta^+}}_{i+1}\le x_{i+1}^\e$. \label{it:mon2}
  \item Let $T= T^{\delta^{-}}:= h^\prime(x_0-\delta)$; then if $ \xie\le x^{\e, T^{\delta^-}}_{i}$, then
    $x_{i+1}^\e\le x^{\e,T^{\delta^-}}_{i+1}$. \label{it:mon3}
  \end{enumerate}
\end{proposition}

\proof
The proof of~\ref{it:mon1} is a straightforward consequence of Proposition~\ref{monof} with
$\phi(t)=\psi(t)= Tt +\e W(t/\e)$ and $\beta = 1/(2\tau_{\e})$.

To prove~\ref{it:mon2}, we apply again Proposition~\ref{monof} with $\phi(t)= h(t) +\e W(t/\e)$ and
$\psi(t)= Tt +\e W(t/\e)$ with $T= T^{\delta^+}:= h^\prime(x_0+\delta)$ and $\beta =
1/(2\tau_{\e})$. By~\eqref{phi-psi}, it holds that
\begin{equation*}
\phi(x_{i+1}^\e) - \phi(x^{\e,T^{\delta^+}}_{i+1})+ \psi(x^{\e,T^{\delta^+}}_{i+1})- \psi(x_{i+1}^\e)\le 2\beta (x^{\e,T^{\delta^+}}_{i}- \xie)
(x^{\e,T^{\delta^+}}_{i+1} - x_{i+1}^\e)\,.
\end{equation*}
Therefore,
\begin{equation*}
h(x_{i+1}^\e)- h(x^{\e,T^{\delta^+}}_{i+1}) + h^\prime(x_0+\delta) (x^{\e,T^{\delta^+}}_{i+1} - x_{i+1}^\e)\le 2\beta (x^{\e,T^{\delta^+}}_{i}- \xie)
(x^{\e,T^{\delta^+}}_{i+1} - x_{i+1}^\e)\,,
\end{equation*}
which implies
 \begin{align*}
& \Biggl(h^\prime(x_0+\delta) - {h(x_{i+1}^\e)- h(x^{\e,T^{\delta^+}}_{i+1})\over  x_{i+1}^\e- x^{\e,T^{\delta^+}}_{i+1} }\Biggr) (x^{\e,T^{\delta^+}}_{i+1} - x_{i+1}^\e)\\
 &\le 2\beta (x^{\e,T^{\delta^+}}_{i}- \xie)
(x^{\e,T^{\delta^+}}_{i+1} - x_{i+1}^\e)\,.
 \end{align*}
 Since $x^{\e,T^{\delta^+}}_{i}, \xie< x_0+\delta$ and $h$ is a convex function,  we get that
\begin{equation*}
  h^\prime(x_0+\delta) - {h(x_{i+1}^\e)- h(x^{\e,T^{\delta^+}}_{i+1})\over  x_{i+1}^\e- x^{\e,T^{\delta^+}}_{i+1} }\ge 0\,,
\end{equation*}
which gives the monotone behaviour $x^{\e,T^{\delta^+}}_{i+1}\le x_{i+1}^\e$. 
\qed


\subsection{Minimising movement for fixed ratio $\e/\tau$}
\label{sec:Minim-movem-fixed}

We consider a time scale $\tau=\tau_\e$ such that $\e/\tau$ converges to $\gamma>0$.  It is not restrictive to suppose
that the ratio between $\e$ and $\tau$ is fixed,
\begin{equation*}
  \tau_\e= {1\over \gamma}\, \e.
\end{equation*} 

We study the linearised energies in~\eqref{F(i)linear} and rescale by $1/\e$, that is,
\begin{equation*}
  {F^T_{\e}(x, \Txie)\over \e} = T\Bigl({x\over \e}\Bigr) + W \Bigl({x\over \e}\Bigr) 
  + {\e\over 2\tau_\e} \Bigl({x- x_{i}^{\e,T}\over \e}\Bigr)^2\,.
\end{equation*}
We denote 
\begin{equation}
  \label{FTgamma}
  F^T_{\gamma}(y, y_i^T):= T y + W (y) + {\gamma\over 2} (y-y^T_{i})^2 \, , 
\end{equation}
where $y:= x/\e$ and $y^T_{i}:= x_{i}^{\e,T}/\e$ for every $i\in\NN$, $i\ge 1$. Note that the minimisers $y^T_{i}$
depend also on $\gamma$. However, we omit this dependence in the notation for simplicity.


\begin{proposition}
  \label{monoy}
  Let $y_0, z_0$ and $T, S$ be fixed with $T\le S$. Let $y_0^T=y_0$, $z_0^S=z_0$.  Let $y_{i}^T$ and $z_{i}^S$ be
  minimisers to $F_{\gamma}^T(y, y^T_{i-1})$ and $F_{\gamma}^S(y, z^S_{i-1})$, respectively, for every $i\in \NN$ with
  $i\ge 1$. If $z_0\le y_0$ then $z_i^{S}\le y_i^{T}$ for every $i$.
\end{proposition}

\proof By Proposition~\ref{monof}, with $\phi(t)= St + W(t)$, $\psi(t)= Tt+W(t)$, $x=z_0$, $x^\prime=y_0$, and
$\beta= \gamma/2$, it follows that
\begin{equation*}
  (S-T) (z^S_1 -y^T_1) \le \gamma (z^S_1 -y^T_1) (z_0-y_0)\,.
\end{equation*}
Therefore, if $T\le S$ and $z_0\le y_0$, this yields $z^S_1 \le y^T_1$. Similarly, we can prove that the inequality
$z_i^{S}\le y_i^{T}$ is satisfied for any $i\ge 2$.  \qed


\begin{theorem}
  \label{tevelo}
  For every $T$, the limit
  \begin{equation}
    \label{fgamma}
    f_{\gamma}(T):=\lim_{i\to \infty} {y_0- y_i^T\over i}
  \end{equation}
  exists and it is independent of $y_0$. Moreover, the function $T\mapsto f_\gamma(T)$ is monotone increasing.
\end{theorem}

\proof 
The existence of the limit is a straightforward consequence of the subadditivity of the sequence $(y_i^T)$. More
precisely, let $h\in \ZZ$ be such that $0\le y_k^T + h\le 1$. Since $ (y_k^T + h)^T_i= y_{k+i}^T + h$, by
Proposition~\ref{monoy}, with $S\equiv T$, 
\begin{equation}
  \label{ineq1}
  y_i^T\le y_{k+i}^T + h\le y_i^T + 1\,.
\end{equation}
Therefore, if we sum up the last inequality in~\eqref{ineq1} with
\begin{equation*}
  -1\le -(y_k^T + h)\le 0\,,
\end{equation*}
we obtain
\begin{equation*}
  y_i^T+ y_k^T-1\le y_{k+i}^T \le y_i^T+ y_k^T+ 1\, ,
\end{equation*}
which implies the almost subadditivity of $(y_i^T)$. We now prove that the limit 
\begin{equation}
  \label{fakete}
  \lim_{i\to\infty} \Bigl({y_i^T\over i}\Bigr) = \inf_{i\in\NN}  \Bigl({y_i^T\over i}\Bigl)
\end{equation}
exists. Let $i= km+n$. Then
\begin{align*}
  \nonumber {y_i^T\over i} & = {y_{km+n}^T\over km+n}
  \le {y_{km}^T+ y_n^T+ 1\over km+n}\\
  \nonumber &\le {k y_{m}^T+ y_n^T+ k\over km+n}= {k y_{m}^T\over km+n}+ { y_n^T\over km+n}+ { k\over km+n}\,. \label{km}
\end{align*}
If we fix $m$ and pass to the limit $k$ tends to $\infty$, we obtain
\begin{equation*}
\lim_{i\to\infty} {y_i^T\over i} \le {y_m^T\over m} + {1\over m}\,.
\end{equation*}
Therefore
\begin{equation*}
\inf_{i\in\NN}{y_i^T\over i}\le \lim_{i\to\infty} {y_i^T\over i} \le \inf_{m\in\NN}{y_m^T\over m} \,,
\end{equation*}
which proves~\eqref{fakete} and the existence of the limit in~\eqref{fgamma}. 

We now prove that the function $T\mapsto f_{\gamma}(T)$ is independent of $y_0$. In fact, we can always rewrite for
$k< i$
\begin{equation*}
  {y_0- y_i^T\over i}= {y_0- y_k^T+ y_k^T- y_i^T\over i} = {y_0- y_k^T\over i} + {y^T_k- y_i^T\over i-k} {i-k\over i}\,.
\end{equation*}
Hence,
\begin{equation*}
  \lim_{i\to \infty} {y_0- y_i^T\over i}= \lim_{i\to \infty} {y^T_k- y_i^T\over i-k}\,.
\end{equation*}
Finally, we remark that the function $T\mapsto f_{\gamma}(T)$ is monotone increasing\ie if $T\le S$ then
$f_{\gamma}(T)\le f_{\gamma}(S)$. By definition~\eqref{fgamma}, 
\begin{equation*}
  f_{\gamma}(T)=\lim_{i\to \infty} {y_0- y_i^T\over i}\,,\qquad f_{\gamma}(S)=\lim_{i\to \infty} {y_0- y_i^{S}\over i}\,;
\end{equation*} 
the monotonicity follows since $- y_i^T\le - y_i^{S}$ by Proposition~\ref{monoy}.  
\qed

\subsection{Characterisation of periodic orbits for the linearised problem}
The definition of $f_\gamma(T)$ reminds that of {\em Poincar\'e   rotation number} in the theory of Dynamical Systems (see, e.g.,~\cite[Chapter 11]{Katok1995a}), which, in our notation, concerns the properties of the orbits of the multifunction
\begin{equation*}
  A^T_\gamma(y)={\rm argmin}\Bigl\{Tz  + W(z) + {1\over 2}\gamma (z-y)^2\Bigr\}.
\end{equation*} 
Note that for $T>1$ this set is a singleton, but for $T\le 1$ in general it is not.  We can nevertheless adapt some
arguments borrowed from Dynamical Systems to prove a characterisation of the values of $T$ for which we have periodic
orbits.  By definition in this case $f_\gamma(T)$ is rational.  The converse also holds true as follows.

\begin{proposition}[Periodic orbits]
  Let $T> 0$, and let $\{y_i^T\}$ be defined as in Proposition~{\rm\ref{monoy}}.  There exists an initial datum
  $y_0=y_0^T$ and integers $p,q$ with $q\neq 0$ such that
  \begin{equation}
    \label{perob}
    y_{kq+i}^T= y^T_i + kp
  \end{equation}
  if and only if $f_\gamma(T)={p\over q}$.
\end{proposition}

\proof 
We only have to prove the existence of $\{y_i^T\}$ satisfying (\ref{perob}) assuming that $f_\gamma(T)=p/q$.
  
We remark that $A=A^T_\gamma$ satisfies \\
$\bullet$ $A$ is monotonically increasing: if $y\le y'$ then $A(y)\le A(y')$; i.e., we have $z\le z'$ for all
$z\in A(y)$ and $z'\in A(y')$ (Proposition~\ref{monof}); \\
$\bullet$ $y\mapsto A(y)$ is (upper) semicontinuous: if $y_n\to y$, $z_n\in A(y_n)$ and $z_n\to z$ then $z\in A(y)$; \\
$\bullet$ $y\mapsto A(y+1)-y$ is $1$-periodic. This last property follows from the $1$-periodicity of $W$, since
\begin{align*}
  A(y+1)&={\rm argmin}\Bigl\{Tz  + W(z) + {1\over 2}\gamma (z-y-1)^2\Bigr\}
  \\
        &=1+{\rm argmin}\Bigl\{T(z-1)  + W(z-1) + {1\over 2}\gamma (z-y)^2\Bigr\}
  \\
        &=1+{\rm argmin}\Bigl\{Tz  + W(z) + {1\over 2}\gamma (z-y)^2\Bigr\}
  \\
        &=1+ A(y).
\end{align*}

Note that the recursive construction of $y^T_i$ translates in $y^T_i\in A(y^T_{i-1})$, and by assumption, we have
\begin{equation*}
  \lim_{n} {y^T_{nq}-y_0\over n}= q f_\gamma(T)= p.
\end{equation*}
Hence, we will examine properties of $y^T_{nq}$, interpreted as the $n$-th iteration of the multifunction $A^q$ (the
$q$-fold composition of $A$) applied to $y_0$.  Note that the multifunction $A^q$ is still increasing and
semicontinuous, and $y\mapsto A^q(y)-y$ is $1$-periodic.

We have to prove that there exists $y_0$ such that
\begin{equation*}
  y_0+p\in A^q(y_0),
\end{equation*}
from which we obtain~\eqref{perob}. Note that we can assume that such $p$ and $q$ are the same as those defining
$f_\gamma(T)$ since this will automatically follow from~\eqref{perob}.

By the monotonicity of $A^q$ we deduce that $A^q(y)$ is a singleton except for a countable number of $y$. We may then
suppose that $A^q(0)$ is a singleton. We denote by $k_0$ the integer part of the unique element of $A^q(0)$, and
consider the multifunction
\begin{equation*}
  G(y):=A^q(y)- k_0.
\end{equation*}
Note that $G$ inherits the properties of $A^q$ and that the unique element of $G(0)$ belongs to $(0,1)$.

We have to show that there exists $y$ such that $(G(y)-y)\cap\ZZ\neq \emptyset$.  We reason by contradiction.  Note
that the graph of $G$ can be extended to a maximal monotone graph $\cal G$ on $\R$, and that if we denote by
${\cal G}(y)$ the corresponding set such that $(y,\bar{y})\in \cal G$ if and only if $\bar{y}\in {\cal G}(y)$, then
${\cal G}(y)$ is a segment (degenerate for almost all $y$) whose endpoints belong to $G(y)$ by semicontinuity.  This
implies that the graph of ${\cal G}(y)-y$ cannot intersect the horizontal lines $\bar{y}\in\ZZ$.  Indeed, suppose that
otherwise there exist $y$ such $0\in {\cal G}(y)-y$ and let $\tilde{y}$ be the minimum of such points in $(0,1)$ (which
exists since the graph of ${\cal G}(y)-y$ is a continuous curve and $0\not\in{\cal G}(0)$.  Then either
${\cal G}(\tilde{y})$ is a singleton, or the segment ${\cal G}(\tilde{y})-\tilde{y}$ has $0$ as the lower endpoint. In
either case, we have $0\in G(\tilde{y})-\tilde{y}$, which contradicts our hypothesis. Similarly, we may show that there
is no $y$ such that $1\in G(y)-y$.  Hence, we have ${\cal G}(y)-y\subset (0,1)$ for all $y$.

By the continuity and periodicity of the graph of ${\cal G}(y)-y$, there exist $\delta>0$ such that
\begin{equation*}
  \delta\le G(y)-y\le 1-\delta\quad  \hbox{for all }y.
\end{equation*}
Let $y_0=0$ and $y_i\in G(y_{i-1})$. For all $n$, from  
\begin{equation*}
  \delta\le y_{i+1}-y_i\le 1-\delta \hbox{ for all }i\in\{0,\ldots, n-1\}
\end{equation*}
we deduce that
\begin{equation*}
  n\delta\le  G(y_n)=A^{nq}(0)-n k_0\le n(1-\delta);
\end{equation*}
that is,
\begin{equation*}
  k_0+ \delta\le {A^{nq}(0)\over n}\le k_0+1- \delta\,.
\end{equation*}
Passing to the limit we finally get
\begin{equation*}
  k_0+ \delta\le q f_\gamma(T)\le k_0+1- \delta,
\end{equation*}
which contradicts the assumption $q f_\gamma(T)\in\ZZ$.
\qed


\section{The limit equation}
\label{sec:limit-equation}

In this section, we show that the limit trajectory $x$ satisfies
\begin{equation}
  \label{lique}
  x'(t)=- \gamma\,f_\gamma(h^\prime(x(t)))
\end{equation}
for almost all $t>0$, with $f_\gamma$ defined in~\eqref{fgamma}. This equation fully characterises $x$ given the
initial datum $x_0$.

\begin{theorem}
  Let $\gamma\in (0, + \infty)$. Let $t_0$ be such that $x^\prime(t_0)$ exists. 
  Then
  \begin{equation*}
    \gamma\,f_\gamma(h^\prime(x(t_0)^-))\le  -x^\prime(t_0)\le  \gamma\,f_\gamma(h^\prime(x(t_0)^+))\,.
  \end{equation*}
\end{theorem}

\proof 
By translating the time variable if necessary we can suppose $t_0=0$ and $x_0= x(0)$. Let $\delta>0$ be such
that $h^\prime(x_0\pm\delta)$ exists. By Proposition~\ref{monotonicity}, \ref{it:mon2}--\ref{it:mon3}, 
\begin{equation*}
  {x^{\e,T^{\delta^+}}_{i} - x_0\over i\tau_\e }\le {{x}^\e_i - x_0\over i\tau_\e }\le {x^{\e,T^{\delta^-}}_{i} - x_0\over  i\tau_\e}\,.
\end{equation*}
The averaged velocity, as in~\eqref{eq:vel}, is given by
\begin{equation*}
  {x_{i}^{\e,T} - x_0\over i \tau_\e }= \gamma {y_i^T- y_0\over i}\,;
\end{equation*}
with definition~\eqref{fgamma} it follows that
\begin{equation*}
  -{\gamma}\, f_\gamma(h^\prime(x_0+\delta))= \lim_{i\to\infty} {x^{\e,T^{\delta^+}}_{i} - x_0\over i\tau_\e }\,,\qquad
  -{\gamma}\, f_\gamma(h^\prime(x_0-\delta))= \lim_{i\to\infty} {x^{\e,T^{\delta^-}}_{i} - x_0\over  i\tau_\e}\,.
\end{equation*}
Therefore, we conclude that
\begin{equation*}
  {\gamma}\, f_\gamma(h^\prime(x_0^-))\le-x^\prime(0) \le {\gamma}\, f_\gamma(h^\prime(x_0^+))
\end{equation*}
as desired.
\qed

The previous result proves that equation~\eqref{lique} fully characterises $x$ when $t\mapsto x(t)$ is strictly
monotone, so that the set of $t$ such that
\begin{equation*}
  \gamma\,f_\gamma(h^\prime(x(t)^-))\neq  \gamma\,f_\gamma(h^\prime(x(t)^+))
\end{equation*}
is of zero (Lebesgue) measure.  
By the monotonicity of $f_\gamma$, if $x$ is not strictly monotone then it is constant, so again~\eqref{lique} is
satisfied. 

We now characterise the \emph{pinning set}, that is, the set of initial data for which $x(t)=x_0$ for all
$t>0$.

\begin{definition}[Pinning threshold]
  For fixed $\gamma>0$, we define the \emph{pinning threshold} at scale $\gamma$ as
  $T_\gamma:=\sup\{T: f_\gamma(T)=0\}$.
\end{definition}

\begin{remark}
  Note that $f_\gamma$ is monotonically increasing, and thus $f_\gamma=0$ on $[0,T_\gamma]$. Hence, for all $x_0$ with
  $|x_0|\le T_\gamma$ the motion is pinned\ie $x(t)=x_0$ for all $t$.
\end{remark}

The following proposition gives a criterion for the computation of the pinning threshold if $Ty+W(y)$ has (at most) a
unique local minimiser in the period. Then it suffices to examine the case where the iteration from that point is
trivial. Note that if $T\ge1$ the function $y\mapsto Ty+W(y)$ is strictly increasing, so that $T>T_\gamma$.

\begin{proposition}[Characterization of the pinning threshold]
  \label{carac}
  Assume that $W'$ has a unique local maximum in $(0,1/2)$. Let $0<T< 1$ and denote by $y_T\in(-1/2,0)$ the unique
  local minimiser of $y\mapsto Ty+W(y)$ in $[-1/2,1/2]$.  Then for every fixed $\gamma>0$, we have $T<T_\gamma$ if and
  only if the function 
  \begin{equation*} 
    \varphi_T(y):= Ty+W(y)+{\gamma\over 2}(y-y_T)^2
  \end{equation*} 
  has a unique global minimum in $y_T$.
\end{proposition}

\proof
Suppose that $y_T$ is the unique global minimiser of $\varphi_T$. Then we can choose as initial datum $y_0=y_T$
in the computation of the velocity in Theorem~\ref{tevelo}, and obtain the trivial orbit $y_k=y_T$. Hence the velocity
is $0$ and consequently $T\le T_\gamma$.  Actually, noting that local minimisers of $\varphi_T$ are a finite set
defined by the identity $T+W'(y)+\gamma(y-y_T)=0$, we have $T<T_\gamma$ by the continuous dependence of these quantities
in $T$.

Conversely, suppose that $y_T$ not be the unique global minimiser of $\varphi_T$. By definition of $T_\gamma$ in order
to show that $T_\gamma\le T$ it suffices to prove that the motion is not pinned for all $T+\delta$ for $\delta>0$.
Then, up to taking such $T+\delta$ in the place of $T$, we may directly suppose that $y_T$ is not a global minimiser of
$\varphi_T$, which is instead a value $y_1\in [y_T-n_1,y_M-n_1]$, where $y_M$ is the unique (local) maximum point of
$W'$ in $(0,1/2)$ and $n_1$ is some positive integer. Now, define the set
\begin{align*}
  I&=\Bigl\{y\in \Bigl(-{1\over 2},{1\over 2}\Bigl): \hbox{ there exists a unique minimiser } \overline y<-{1\over 2}\  \hbox {of } \\
 &\qquad\qquad\qquad\qquad\qquad\qquad\quad w\mapsto Tw+W(w)+{\gamma\over 2}(w-y)^2\Bigr\},
\end{align*}
which is the set of initial data for which the first iteration moves to ``another well''.  By continuity, there exists
$\delta>0$ such that $[y_T-\delta,y_T+\delta]\subset I$. Then
$y_{N}\in I-k_1$ after a finite number of iterations $N$ independent of $y_1$ (for a finer estimate of $N$ in the
piecewise-quadratic case we refer to Section~\ref{sec:An-exampl-piec}), and we can proceed by induction.  This gives
the positiveness of the velocity and $T\ge T_\gamma$.
\qed

\begin{remark}[Asymptotic behaviour at the pinning threshold]
  \label{asybe}
  If the hypotheses of Proposition~\ref{carac} are satisfied and $W$ is $C^2$ at local minimisers of $Ty+W(y)$ with
  strictly positive second derivative then for all $\gamma>0$
  \begin{equation*}
    f_\gamma(T)\sim {1\over \log(T-T_\gamma)}
  \end{equation*}
  as $T\to T_\gamma^+$. This will be shown for piecewise quadratic energies $W$ in detail in the next section.
\end{remark}

\begin{proposition}[Extreme minimizing movements] 
  \label{prop:exmmm}
  We have
  \begin{equation}
    \lim_{\gamma\to 0} \gamma f_\gamma(z) = z,\qquad
    \lim_{\gamma\to +\infty} \gamma f_\gamma(z) = g_\infty(z),
  \end{equation}
  where $g_\infty$ is given by
  \begin{equation}
    \label{ginfty}
    g_\infty(z)=
    \begin{cases}\displaystyle
      \Biggl(\int_0^1{1\over z+W(s)}ds\Biggl)^{-1} &\text{ if }\displaystyle{1\over z+W(s)} \text{ is integrable,}\cr\cr
      0 & \text{otherwise}.
    \end{cases}
  \end{equation} 
  Moreover, 
  \begin{equation}
    \lim_{\gamma\to+\infty} T_\gamma= \sup_\gamma T_\gamma=T_\infty,
  \end{equation}
  where $[-T_\infty,T_\infty]=\{T\in\R: g_\infty(T)=0\}$.
\end{proposition}

\proof 
Assume that $h$ and $W$ are $C^2$-functions. The convergence as $\gamma\to 0$ follows from the observation that the orbit $x^\e_k$ satisfies
\begin{equation*}
  {x^\e_k-x^\e_{k-1}\over\tau}= -h^\prime(x^\e_k) +O\Bigl({\e\over\tau}\Bigr)
  = -x^\e_k +o(1)
\end{equation*}
as $\gamma\to 0$. Conversely, the convergence as $\gamma\to +\infty$ follows by noting that as $\gamma\to +\infty$ the
orbits $x^\e(t)$, defined as in (\ref{seqMM}), are close to the corresponding solution of the gradient flow
\begin{equation*}
  x_\e'=-h'(x_\e)-W^\prime\Bigl({x_\e\over\e}\Bigl),
\end{equation*}
whose limit satisfies $x'=-g_\infty(h'(x))$.
\qed

\begin{remark}
  \label{extra}
  By Theorem 8.1 in~\cite{Braides2014a}, the equations
  \begin{equation*}
    x'=-h'(x) \text{ and } x'=-g_\infty(h'(x))
  \end{equation*}
  describe the minimising movements in the cases $\e\ll\tau$ and $\tau\ll\e$, respectively.  The previous proposition
  shows that the same extreme minimising movements are obtained by keeping the ratio $\gamma=\e/\tau$ fixed and then
  let it tend to $0$ and $+\infty$, respectively.
\end{remark}


\section{An example: the piecewise-quadratic case}
\label{sec:An-exampl-piec}

In this section, we provide an example of oscillating potential and calculate explicitly the corresponding pinning
threshold $T_\gamma$. More precisely, we consider the piecewise quadratic energy
\begin{equation*}
  W(y):= \min_{k\in \ZZ} (y-k)^2\,.
\end{equation*}

Besides giving an illustrative example, we deduce the asymptotic behaviour at the pinning threshold, which depends only
on the non-degeneracy of the second derivative at local minima. In this way we deduce the asymptotic behaviour in
the general case as in Remark~\ref{asybe}.

For this choice of $W$, for $T_\infty$ as in Proposition~\ref{prop:exmmm}, it holds that
\begin{equation*}
  T_\infty= 1,
\end{equation*}
and~\eqref{ginfty} becomes
\begin{equation*}
  g_\infty(z)={1\over \log\bigl({z-1\over z}\bigr)} \qquad\hbox{ for } z> 1.
\end{equation*}

The function $W$ is $1$-periodic and piecewise quadratic with $\Vert W^\prime\Vert_\infty=1$.  For simplicity, we fix
also $h(x)= x^2/2$. Let $T\in (0,1)$ and $y_0\in [0, 1/2)$. The minimum of the function
$F^T_{\gamma}(y, y_0)= T y + W (y) + {\gamma\over 2} (y-y_0)^2$, given in~\eqref{FTgamma}, on the interval
$[-1/2 +k, 1/2 +k]$ can be attained at the boundary or at the interior of this interval; it is given by
\begin{equation}
  \label{formulayk}
  y_{1,k}= {-T+2k\over 2+\gamma} + {\gamma\over 2+\gamma} y_0
\end{equation} 
($y_{1,k}$ also depends on $T$ but we suppress this in the notation).

The global minimiser $y^T_1$ to $F^T_{\gamma}(y, y_0)$ can get stuck in the same well of $y_0$, that is,
$y^T_1= y_{1,0}$. Otherwise it can move into the next well, corresponding to $k=-1$, that is, $y^T_1= y_{1,-1}$.

Any single well of $W(y)$ is denoted by $W(y;k):= (y-k)^2$ for every $y\in [-1/2+k, 1/2+k]$ with $k\in \ZZ$.  We define
\begin{align}
  \psi (y) &:= \Bigl(Ty_{1,-1}+ W(y_{1,-1};-1) +  {\gamma\over 2} (y_{1,-1}- y)^2 \Bigr) \notag\\
 & {}\qquad- \Bigl( Ty_{1,0}+ W(y_{1,0};0) +  {\gamma\over 2} (y_{1,0}- y)^2
  \Bigr) \,.  \label{defpsi}
\end{align}
To establish if $y_1^T$ gets stuck or moves, we have to study the sign of $\psi (y_0)$, since it is the difference
between the minimum value of the two wells. Therefore, if $\psi (y_0) < 0$, then the minimiser satisfies
$y^T_1= y_{1,-1}$; $\psi (y_0) \ge 0$ implies that $y^T_1= y_{1,0}$. In particular, from the sign of $\psi (y_0)$ we
expect to derive the pinning threshold $T_\gamma$.


\begin{proposition}
  \label{quaqua}
  Let $\gamma\in (0, +\infty)$ and let $T\in(0,1)$. Then there exist
  \begin{equation*}
    T_\gamma := {\gamma\over  (2+\gamma)} \quad   \hbox{and}   \quad  \delta_T:= \Bigl({2+\gamma\over 2\gamma}\Bigr)\, (T- T_\gamma),
  \end{equation*}
  such that the following holds.
  \begin{enumerate}
  \item \label{it:quaqua1}
    For every $T>T_\gamma$ we have that, the following possibilities exist.
    \begin{enumerate}
    \item \label{it:quaqua1a} If $y_0\in [0, (-T/2)+\delta_T)$ then $y_1^T= y_{1,-1}$, where the latter is defined
      in~\eqref{formulayk}. Moreover, if
      \begin{equation*}
        y_1^T +{T\over 2}+1=  {\gamma\over 1+\gamma} \Bigl(y_0 + {T\over 2}+1\Big) < -{T\over 2} +\delta_T\,,
      \end{equation*}
      then the successive minimiser is given by $y^T_2= y_{2,-2}$ and so on.
    \item \label{it:quaqua1b} If $y_0\in [(-T/2)+\delta_T, 1/2)$
      then there exists $h\in \NN$  given by
      \begin{equation}\label{def-h}
        h= \Bigg\lfloor{\log \Bigl( ({2+\gamma\over \gamma}) {T- T_\gamma \over T+ 1}\Bigr) 
          \over \log \Bigl({\gamma\over 2+\gamma}\Bigr)}\Bigg\rfloor + 1 
      \end{equation}
      such that $y_1^T= y_{1,0}, \ldots,y_{h}^T= y_{h-1,0} \ge (-T/2) +\delta_T$, $y_{h,0} < (-T/2) +\delta_T$ and
      $y_{h+1}^T= y_{h+1,-1}$. Similarly, if $k\in\NN$ exists such that for some $p\in\NN$
      \begin{equation*}
        y_k^T +{T\over 2}+p \ge  -{T\over 2} +\delta_T \,,
      \end{equation*} 
      then there exists $h$ as in~\eqref{def-h} such that
      if we take as initial data $z_0:= y^T_k +(T/2)+p$ then we get a new sequence of minimisers such that
      $z^T_1=z_{1,0}\,,\ldots ,z^T_{h-1}= z_{h-1,0}\ge -T +\delta_T$, $z_{h}^T= z_{h,0} < (-T/2) +\delta_T$ and
      $z_{h+1}^T= z_{h+1,-1}$.
    \end{enumerate}
  \item \label{it:quaqua2} For every 
    \begin{equation*}
      T\le  T_\gamma\ \hbox{and}\ y_0\ge 0
    \end{equation*}  the motion is pinned.
  \end{enumerate}
\end{proposition}

\proof We first derive $T_\gamma$ using the criterion given in Proposition~\ref{carac}. More precisely, since the
unique local minimiser of $y\mapsto Ty+W(y)$ in $[-1/2,1/2]$ is $y_T=-T/2$, by~\eqref{formulayk} and~\eqref{defpsi},
with $y_0=y_T$ we have that
\begin{equation*}
  y_{1,-1}= -{T\over 2} -{2\over 2+\gamma} \,,\qquad 
  y_{1,0}= -{T\over 2} \,.
\end{equation*} 
and
\begin{equation*}
  \psi (y_T)= {2\over 2+\gamma} \Bigl(-T + {\gamma\over 2} + \gamma y_T\Bigr)=
  -T +{\gamma\over 2 +\gamma}\,.
\end{equation*}
Therefore, we find
\begin{equation*}
  T_\gamma= {\gamma\over 2 +\gamma}\,.
\end{equation*}
Moreover, for fixed $y_0= (-T/2) +\delta$, one has
\begin{equation*}
  \psi\Bigl(-{T\over 2} +\delta\Bigr)\le 0 \quad \hbox{if and only if} \quad  
  \delta\le {2+\gamma\over 2\gamma} \Bigl(T-{\gamma\over 2+\gamma}\Bigr)\,.
\end{equation*}
Therefore,  if we define 
\begin{equation}
  \label{deltaT}
  \delta_T:=  {2+\gamma\over 2\gamma} \Bigl( T- T_\gamma\Bigr)
\end{equation}
we are ready to prove the statements~\ref{it:quaqua1a} and~\ref{it:quaqua1b} of the proposition.

We now give the proof of~\ref{it:quaqua1a}. For every $T > T_\gamma$ there exists $\delta_T>0$, given
by~\eqref{deltaT}, such that for every $y_0\in [0, (-T/2)+\delta_T)$ we have that $\psi (y_0) <0$, that is,
$ y^T_1= y_{1,-1}$.

Reasoning as above, we observe that if $y_{1,-1} +(T/2)+1< (-T/2) +\delta_T$, then $y^T_2= y_{2,-2}$ and we can iterate
until this condition is satisfied. The case $y_{1,-1} +(T/2)+1\ge (-T/2) +\delta_T$ is addressed in
point~\ref{it:quaqua1b}.

Next, we give the proof of claim~\ref{it:quaqua1b}. If $y_0\in [(-T/2)+\delta_T, 1/2)$ then $\psi (y_0)\ge0$ and the
minimisers can be calculated recursively, by~\eqref{formulayk}, in the following way.
\begin{align*}
 y^T_1&=y_{1,0} = {\gamma\over 2+\gamma} \, y_0 - {T\over 2+\gamma}\,,  \\
 y^T_2 &=y_{2,0} = {\gamma\over 2+\gamma} \, y^T_1 - {T\over 2+\gamma}   \\
 & \qquad = \Bigl({\gamma\over 2+\gamma}\Bigr)^2 \, y_0 - {T\over 2+\gamma} \Bigl( 1 + {\gamma\over 2+\gamma}\Bigr)\,,\\
  y^T_h &= y_{h,0} 
 = \Bigl({\gamma\over 2+\gamma}\Bigr)^h y_0 - {T\over 2+\gamma}  \sum_{n=0}^{h-1} \Bigl( {\gamma\over 2+\gamma}\Bigr)^n\,.
\end{align*}
Therefore, we may rewrite
\begin{equation}
  \label{recursive-yh}
  y^T_{h} = \Bigl({\gamma\over 2+\gamma}\Bigr)^h \Bigl(y_0 + {T\over 2}\Bigr) -{T\over 2}\,.
\end{equation}
Since $(y_0 + (T/2)) < {1\over 2} +(T/2)$, then we may assume that
\begin{equation*}
\Bigl({\gamma\over 2+\gamma}\Bigr)^h \Bigl(y_0 + {T\over 2}\Bigr) -{T\over 2} <
\Bigl({\gamma\over 2+\gamma}\Bigr)^h \Bigl({T+1\over 2} \Bigr) -{T\over 2} <  -{T\over 2}+\delta_T\,.
\end{equation*}
Therefore
\begin{equation*}
  \Bigl({\gamma\over 2+\gamma}\Bigr)^h < {2\delta_T\over T + 1};
\end{equation*} 
that is,
\begin{align*}
  h &> {\log \Bigl( {2\delta_T\over T+ 1}\Bigr) \over \log \Bigl({\gamma\over 2+\gamma}\Bigr)} =
      {\log \Bigl(  ({2+\gamma\over \gamma}) {T- T_\gamma \over T+ 1}\Bigr) \over \log \Bigl({\gamma\over 2+\gamma}\Bigr)}\,.
\end{align*}
Similarly, if $y^T_1 +(T/2)+1 \ge -(T/2) +\delta_T$ then we may reason as above, by taking as initial datum
$z_0:= y^T_1 +(T/2)+1$. Therefore, we get a new sequence of minimisers given by
\begin{align*}
z^T_1&=z_{1,0} = {\gamma\over 2+\gamma} \, z_0 - {T\over 2+\gamma} \,, \\
z^T_2 &=z_{2,0} = {\gamma\over 2+\gamma} \, z_1 - {T\over 2+\gamma}   \\
& \qquad = \Bigl({\gamma\over 2+\gamma}\Bigr)^2 \, z_0 - {T\over 1+\gamma} \Bigl( 1 + {\gamma\over 2+\gamma}\Bigr) \,,\\
 & \cdots  \\
z^T_h &= z_{h,0} 
= \Bigl({\gamma\over 2+\gamma}\Bigr)^h \Bigl(z_0 + {T\over 2}\Bigr) -{T\over 2}\,,
\end{align*}
such that $z^T_1=z_{1,0}\,,\ldots ,z^T_{h-1}= z_{h-1,0}\ge (-T/2) +\delta_T$, $z_{h}^T= z_{h,0} < (-T/2) +\delta_T$
and $z_{h+1}^T= z_{h+1,-1}$.  More generally, if $k\in\NN$ exists such that, for some $p\in\NN$,
$y^T_k +(T/2)+p \ge (-T/2) +\delta_T$, then we may repeat the procedure above by assuming $z_0:= y^T_k +(T/2)+p$.
 

We now turn to case~\ref{it:quaqua2} and give the proof. If $T\le T_\gamma$, since $y_0\ge 0$, then $\psi(y_0)>0$, that
is, $y_1^T= y_{1,0}$. Moreover, by~\eqref{recursive-yh}, we have that
\begin{equation*}
y_h^T= \Bigl({\gamma\over 2+\gamma}\Bigr)^h \Bigl(y_0 + {T\over 2}\Bigr) -{T\over 2} > -{T\over 2} +\delta_T
\end{equation*}
for every $h\in\NN$. Therefore $\lim_{h\to \infty} y_h^T= (-T/2)$, that is, for every $T\le T_\gamma$ the motion is
pinned. \qed


\begin{remark}[Behaviour at the pinning threshold]
  From Proposition~\ref{quaqua}, case~\ref{it:quaqua1b} we deduce that
  \begin{equation*}
    f_\gamma (T)\sim  {\log \bigl({\gamma\over 2+\gamma}\bigr)\over\log \Bigl(({2+\gamma\over \gamma}) {T- T_\gamma \over T+ 1}\Bigr)  };
  \end{equation*}
  that is, for $\gamma>0$ fixed,
  \begin{equation*}
    f_\gamma (T)\sim  {1\over|\log (T- T_\gamma) | }
  \end{equation*}
  as $T\to T_\gamma^+$. Note in particular that $f_\gamma$ is not Lipschitz for $T\to T_\gamma^+$.
\end{remark}


\paragraph{Acknowledgements} NA and JZ gratefully acknowledge funding by the Marie Curie Actions: Intra-European
Fellowship for Career Development (IEF2012, FP7-People) under REA grant agreement n${}^\circ$ 326044. JZ was partially
supported by the UK's Engineering and Physical Sciences Research Council Grant EP/K027743/1 and the Leverhulme Trust
(RPG-2013-261).


\begin{thebibliography}{10}

\bibitem{Abeyaratne1996a}
R.~Abeyaratne, C.~Chu, and R.~D. James.
\newblock Kinetics of materials with wiggly energies: theory and application to
  the evolution of twinning microstructures in a {C}u-{A}l-{N}i shape memory
  alloy.
\newblock {\em Philos. Mag. A}, 73(2):457--497, February 1996.

\bibitem{Abeyaratne1999a}
R.~Abeyaratne and S.~Vedantam.
\newblock Propagation of a front by kink motion---from a discrete model to a
  continuum model.
\newblock In {\em Variations of domain and free-boundary problems in solid
  mechanics (Paris, 1997)}, volume~66 of {\em Solid Mech. Appl.}, pages 77--84.
  Kluwer Acad. Publ., Dordrecht, 1999.

\bibitem{Ambrosio2008b}
L.~Ambrosio, N.~Gigli, and G.~Savar{\'e}.
\newblock {\em Gradient flows in metric spaces and in the space of probability
  measures}.
\newblock Lectures in Mathematics ETH Z\"urich. Birkh\"auser Verlag, Basel,
  second edition, 2008.

\bibitem{Braides2014a}
A.~Braides.
\newblock {\em Local Minimization, Variational Evolution and
  {$\Gamma$}-con\-vergence}, volume 2094 of {\em Lecture Notes in Mathematics}.
\newblock Springer, Cham, 2014.

\bibitem{Braides2008a}
A.~Braides and L.~Truskinovsky.
\newblock Asymptotic expansions by $\Gamma$-convergence.
\newblock {\em Continuum Mech. Thermodyn.}, 20:21--62, 2008.

\bibitem{De-Giorgi1993a}
E.~De~Giorgi.
\newblock New problems on minimizing movements.
\newblock In {\em Boundary value problems for partial differential equations
  and applications}, volume~29 of {\em RMA Res. Notes Appl. Math.}, pages
  81--98. Masson, Paris, 1993.

\bibitem{Dupuis2012a}
P.~Dupuis and K.~Spiliopoulos.
\newblock Large deviations for multiscale diffusion via weak convergence
  methods.
\newblock {\em Stochastic Process. Appl.}, 122(4):1947--1987, 2012.

\bibitem{Ibrahim2010a}
H.~Ibrahim and R.~Monneau.
\newblock On the rate of convergence in periodic homogenization of scalar
  first-order ordinary differential equations.
\newblock {\em SIAM J. Math. Anal.}, 42(5):2155--2176, 2010.

\bibitem{Katok1995a}
A.~Katok and B.~Hasselblatt.
\newblock {\em Introduction to the modern theory of dynamical systems},
  volume~54 of {\em Encyclopedia of Mathematics and its Applications}.
\newblock Cambridge University Press, Cambridge, 1995.
\newblock With a supplementary chapter by Katok and Leonardo Mendoza.

\bibitem{Menon2002a}
G.~Menon.
\newblock Gradient systems with wiggly energies and related averaging problems.
\newblock {\em Arch. Ration. Mech. Anal.}, 162(3):193--246, 2002.

\end{thebibliography}

\def\cprime{$'$} \def\cprime{$'$} \def\cprime{$'$}
  \def\polhk#1{\setbox0=\hbox{#1}{\ooalign{\hidewidth
  \lower1.5ex\hbox{`}\hidewidth\crcr\unhbox0}}} \def\cprime{$'$}
  \def\cprime{$'$}

\end{document}